\documentclass[english]{amsart}

\usepackage{amsmath, amssymb}
\usepackage{etoolbox} 
\usepackage{mathtools}
\usepackage{interval} 
\usepackage{hyperref}
\intervalconfig{soft open fences}

\newcommand\numberthis{\addtocounter{equation}{1}\tag{\theequation}}

\newcommand{\Z}{\ensuremath{\mathbb{Z}}}

\let\Re\relax
\let\Im\relax
\DeclareMathOperator{\Re}{Re}
\DeclareMathOperator{\Im}{Im}

\DeclareMathOperator{\SL}{SL}
\DeclareMathOperator{\GL}{GL}

\newtheorem{theorem}{Theorem}
\newtheorem{remark}[theorem]{Remark}
\newtheorem{lemma}[theorem]{Lemma}

\makeatletter
\let\@@pmod\pmod
\DeclareRobustCommand{\pmod}{\@ifstar\@pmods\@@pmod}
\def\@pmods#1{\mkern4mu({\operator@font mod}\mkern 6mu#1)}
\makeatother

\newcommand*{\conjugate}[1]{\overline{#1}}

\newcommand*{\placeholder}{\makebox[1ex]{\textbf{$\cdot$}}}

\DeclarePairedDelimiterX{\abs}[1]{\lvert}{\rvert}{
	\ifblank{#1}{\placeholder}{#1}%
}

\DeclarePairedDelimiterX{\set}[1]{\{}{\}}{
	\,#1\,
}

\let\oldset\set
\def\set{\futurelet\testchar\MaybeOptArgSet}
\def\MaybeOptArgSet{\ifx[\testchar \let\next\OptArgSet
\else \let\next\NoOptArgSet \fi \next}
\def\OptArgSet[#1]#2{\oldset[#1]{#2}}
\def\NoOptArgSet#1{\OptArgSet[\big]{#1}}

\makeatletter
\newcommand{\Set}{%
             \@ifstar
                  {\oldset*}%
                  {\set}%
}
\makeatother

\let\originalleft\left
\let\originalright\right
\renewcommand{\left}{\mathopen{}\mathclose\bgroup\originalleft}
\renewcommand{\right}{\aftergroup\egroup\originalright}

\begin{document}

\title{Twisted Moments of $\GL(3) \times \GL(2)$ $L$-functions}
\author{Jakob Streipel}

\address{Department of Mathematics and Statistics, Washington State University, Pullman, WA 99164, USA}
\email{jakob.streipel@wsu.edu}

\begin{abstract}
	We compute an asymptotic formula for the twisted moment of $\GL(3) \times \GL(2)$ $L$-functions and their derivatives.
	As an application we prove that symmetric-square lifts of $\GL(2)$ Maass forms are uniquely determined by the central values of the derivatives of $\GL(3) \times \GL(2)$ $L$-functions.
\end{abstract}

\maketitle



\section{Introduction}

The method of moments has broad applications to the study of $L$-functions and the forms attached to them.
In her groundbreaking paper \cite{li2011}, Li uses bounds for moments of $\GL(3) \times \GL(2)$ $L$-functions to deduce the first subconvexity bound for degree three $L$-functions.
Using asymptotic formulas for these moments one can derive nonvanishing results.
For example \cite{li2009, sun2012, chenyan2018} all use asymptotics for moments of $\GL(3) \times \GL(2)$ $L$-functions or their derivative to deduce simultaneous nonvanishing results at the central point of these objects.
In a series of papers Liu (\cite{Liu2010, Liu2011}) and later Sun (\cite{sun2014,sun2015}) use asymptotics for these moments twisted by Fourier coefficients at primes to show that self-dual Maass cusp forms for $\GL(3)$ are uniquely determined by central values of $\GL(3) \times \GL(2)$ $L$-functions.

In this paper we compute asymptotics for twisted moments of $\GL(3) \times \GL(2)$ $L$-functions (\autoref{thm:maintheorem}) and their derivatives (\autoref{thm:derivativemoment}).
As an application of this we prove an analogue of Liu's results in the spectral aspect, i.e.\ we show that symmetric-square lifts of $\GL(2)$ Hecke--Maass forms are uniquely determined by central values of the derivatives of $\GL(3) \times \GL(2)$ $L$-functions (\autoref{thm:determination}).

%
%
%

Let $f$ be a Hecke--Maass form of type $(\nu_1, \nu_2)$ for $\SL(3, \Z)$ with Fourier coefficients $A(m, n)$, normalized so that the first Fourier coefficient is $A(1, 1) = 1$.
We define the $L$-function
\[
	L(s, f) = \sum_{m = 1}^\infty \frac{A(m, 1)}{m^s}
\]
for $\Re(s) > 1$.
This has analytic continuation and its completed $L$-function $\Lambda(s, f)$ satisfies $\Lambda(s, f) = \Lambda(1 - s, \tilde f)$, where $\tilde f$ denotes the dual form of $f$ of type $(\nu_2, \nu_1)$ and Fourier coefficients $A(n, m)$.
Let $\Set{u_j}$ be an orthonormal basis of Hecke--Maass forms for $\SL(2, \Z)$.
We define the Rankin--Selberg $L$-function
\[
	L(s, f \times u_j) = \sum_{m \geq 1} \sum_{n \geq 1} \frac{\lambda_j(n) A(n, m)}{(m^2 n)^s}
\]
for $\Re(s) > 1$, where $\lambda_j(n)$ are the normalized Fourier coefficients of $u_j$ (see \nameref{sec:prelim} for details).

Our main result is this:

\begin{theorem}\label{thm:maintheorem}
	Let $f$ be a Hecke--Maass form for $\SL(3, \Z)$ and let $\Set{u_j}$ be an orthonormal basis of even Hecke--Maass forms for $\SL(2, \Z)$ with Laplacian eigenvalues $\frac{1}{4} + t_j^2$, $t_j \geq 0$ and normalized Fourier coefficients $\lambda_j(n)$.
	Fix a prime $p \ll T^{1 - \varepsilon}$.

	Then for any $\varepsilon > 0$ and $T$ large with $T^{\frac{3}{8} + \varepsilon} < M \leq T^{1 - \varepsilon}$, we have
	\begin{align*}
		&\sideset{}{'}\sum_{j} k(t_j) \omega_j \lambda_j(p) L\Bigl( \frac{1}{2}, f \times u_j \Bigr ) + \frac{1}{4 \pi} \int_{-\infty}^\infty k(t) \omega(t) \conjugate{\eta}\Bigl ( p, \frac{1}{2} + i t \Bigr ) \abs[\Big]{L\Bigl ( \frac{1}{2} - i t, f \Bigr ) }^2 \, d t \\
		&\qquad = \frac{L(1, \tilde f) \bigl(A(p, 1) p - 1\bigr) + L(1, f) \bigl(A(1, p) p - 1\bigr)}{p^{\frac{3}{2}} \pi} \int_0^\infty k(t) \tanh(\pi t) t \, d t + {} \\
		&\qquad\qquad{} + O(M^{-3} T^{\frac{5}{2} + \varepsilon} p^{1 + \varepsilon} + M^{-1} T^{\frac{3}{2} + \varepsilon} p^{\varepsilon} + M T^{\frac{1}{7} + \varepsilon} p^{\varepsilon}),
	\end{align*}
	where $'$ means summing over the orthonormal basis of even Hecke--Maass forms, $k(t) = e^{-\frac{(t - T)^2}{M^2}} + e^{-\frac{(t + T)^2}{M^2}}$, and the weights $\omega_j$ and $\omega(t)$ as well as $\eta(n, s)$ are defined in the \nameref{sec:prelim}.
\end{theorem}

\begin{remark}
	The integral in the main term is of size
	\[
		\int_0^\infty k(t) \tanh(\pi t) t \, d t \asymp T M.
	\]
\end{remark}

\begin{remark}
	Without the twist $\lambda_j(p)$ (i.e., corresponding to $p = 1$), the shape of the asymptotic formula is essentially the same, with main term
	\[
		\frac{L(1, \tilde f) + L(1, f)}{\pi} \int_0^\infty k(t) \tanh(\pi t) t \, d t.
	\]
\end{remark}

If $f$ is self-dual this becomes
\begin{align*}
	&\sideset{}{'}\sum_{j} k(t_j) \omega_j \lambda_j(p) L\Bigl( \frac{1}{2}, f \times u_j \Bigr ) + \frac{1}{4 \pi} \int_{-\infty}^\infty k(t) \omega(t) \conjugate{\eta}\Bigl ( p, \frac{1}{2} + i t \Bigr ) \abs[\Big]{L\Bigl ( \frac{1}{2} - i t, f \Bigr ) }^2 \, d t \\
	&\qquad= \frac{2 L(1, f) \bigl (A(p, 1) p - 1 \bigr)}{p^{\frac{3}{2}} \pi} \int_0^\infty k(t) \tanh(\pi t) t \, d t + {} \\
	&\qquad\qquad{} + O(M^{-3} T^{\frac{5}{2} + \varepsilon} p^{1 + \varepsilon} + M^{-1} T^{\frac{3}{2} + \varepsilon} p^{\varepsilon} + M T^{\frac{1}{7} + \varepsilon} p^{\varepsilon}).
\end{align*}
\begin{remark}
	Note that the analysis of the last error term, coming from the diagonal, can be refined in the case where $f$ is self-dual (using stronger bounds on the Ramanujan conjecture in this case), however the improvement in the error term doesn't affect our application of it below.
\end{remark}

Our approach in proving this, contained in Sections~\ref{sec:setup}--\ref{sec:Kbessel}, is based on \cite{li2011}.

Using the same technology we further compute, in Section~\ref{sec:derivativemoment}, asymptotics for the twisted first moment of the derivative at the central point:

\begin{theorem}\label{thm:derivativemoment}
	Let $f$ be a Hecke--Maass form for $\SL(3, \Z)$ and let $\Set{u_j}$ be an orthonormal basis of odd Hecke--Maass forms for $\SL(2, \Z)$ with Laplacian eigenvalues $\frac{1}{4} + t_j^2$, $t_j \geq 0$ and normalized Fourier coefficients $\lambda_j(n)$.
	Fix a prime $p \ll T^{1 - \varepsilon}$.
	Then for any $\varepsilon > 0$ and $T$ large with $T^{\frac{3}{8} + \varepsilon} < M \leq T^{1 - \varepsilon}$, we have
	\begin{align*}
		&\sideset{}{^*}\sum_j k(t_j) \omega_j \lambda_j(p) L'\Bigl( \frac{1}{2}, f \times u_j \Bigr) \\
		&= \frac{3 L(1, \tilde f) \bigl( A(p, 1) p - 1\bigr) + 3 L(1, f) \bigl( A(1, p) p - 1\bigr)}{2 p^{\frac{3}{2}} \pi} \int_{-\infty}^\infty k(t) \tanh(\pi t) t \log\abs{t} \, d t + {} \\
		&\qquad{} + \frac{K \bigl( A(p, 1) p - 1\bigr) + \tilde K \bigl( A(1, p) p - 1\bigr)}{2 p^{\frac{3}{2}} \pi} \int_{-\infty}^{\infty} k(t) \tanh(\pi t) t \, d t + {} \\
		&\qquad\qquad{} + O(M^{-3} T^{\frac{5}{2} + \varepsilon} p^{1 + \varepsilon} + M^{-1} T^{\frac{3}{2} + \varepsilon} p^{\varepsilon} + M T^{\frac{1}{7} + \varepsilon} p^{\varepsilon})
	\end{align*}
	where $^*$ means summing over the orthonormal basis of odd Hecke--Maass forms, $K = 2 L'(1, \tilde f) - 3 L(1, \tilde f) \log(2 \pi) - L(1, \tilde f) \log p$ and $\tilde K = 2 L'(1, f) - 3 L(1, f) \log(2 \pi) - L(1, f) \log p$ are constants, $k(t) = e^{-\frac{(t - T)^2}{M^2}} + e^{-\frac{(t + T)^2}{M^2}}$, and the weight $\omega_j$ is defined in the \nameref{sec:prelim}.
\end{theorem}

\begin{remark}
	The main term is of size
	\[
		\int_{-\infty}^\infty k(t) \tanh(\pi t) t \log\abs{t} \, d t \asymp T M \log T.
	\]
\end{remark}

\begin{remark}
	Without the twist $\lambda_j(p)$, the asymptotic is again similar, with main term
	\[
		\frac{3 L(1, \tilde f) + 3 L(1, f)}{2 \pi} \int_{-\infty}^\infty k(t) \tanh(\pi t) t \log\abs{t} \, d t.
	\]
	This agrees with \cite[Theorem~1.1]{chenyan2018} when taking the average as long as possible (i.e., $M$ close to $T$), but also allows averaging over shorter intervals.
\end{remark}

For self-dual $f$ this becomes
\begin{align*}
	&\sideset{}{^*}\sum_j k(t_j) \omega_j \lambda_j(p) L'\Bigl( \frac{1}{2}, f \times u_j \Bigr) \\
	&= \frac{3 L(1, f) \bigl( A(p, 1) p - 1\bigr)}{p^{\frac{3}{2}} \pi} \int_{-\infty}^\infty k(t) \tanh(\pi t) t \log\abs{t} \, d t + {} \\
	&\qquad{} + \frac{K \bigl( A(p, 1) p - 1\bigr)}{p^{\frac{3}{2}} \pi} \int_{-\infty}^{\infty} k(t) \tanh(\pi t) t \, d t + {} \\
	&\qquad\qquad{} + O(M^{-3} T^{\frac{5}{2} + \varepsilon} p^{1 + \varepsilon} + M^{-1} T^{\frac{3}{2} + \varepsilon} p^{\varepsilon} + M T^{\frac{1}{7} + \varepsilon} p^{\varepsilon}).
\end{align*}


An application of this is to show that symmetric-square lifts of $\GL(2)$ Maass forms are uniquely determined by the central values $L'(\frac{1}{2}, f \times u_j)$.

\begin{theorem}\label{thm:determination}
	Let $f$ and $g$ be symmetric-square lifts of $\GL(2)$ Maass forms.
	If
	\[
		L'\Bigl( \frac{1}{2}, f \times u_j \Bigr) = c L'\Bigl( \frac{1}{2}, g \times u_j \Bigr)
	\]
	for some constant $c \neq 0$ and all odd $u_j$, then $f = g$.
\end{theorem}

This is an analogue in the spectral aspect of Liu's results \cite{Liu2010,Liu2011} in the weight and level aspects.


We quickly sketch the proof of this.
Let $A_f(m, n)$ and $A_g(m, n)$ denote the Fourier coefficients of $f$ and $g$ respectively.
For each fixed prime $p$, the assumption of $L'(\frac{1}{2}, f \times u_j) = c L'(\frac{1}{2}, g \times u_j)$ applied to \autoref{thm:derivativemoment} tells us
\[
	\frac{3 L(1, f)}{\pi} \Bigl( \frac{A_f(p, 1)}{p^{\frac{1}{2}}} - \frac{1}{p^{\frac{3}{2}}} \Bigr) T M \log T = \frac{3 c L(1, g)}{\pi} \Bigl( \frac{A_g(p, 1)}{p^{\frac{1}{2}}} - \frac{1}{p^{\frac{3}{2}}} \Bigr) T M \log T + O_p(T M).
\]
Moreover computing the same kind of moment as in \autoref{thm:derivativemoment} without the twist by $\lambda_j(p)$ we have
\[
	L(1, f) T  M \log T = c L(1, g) T M \log T + O_p(T M).
\]
By taking $T \to \infty$ these two together therefore imply $A_f(p, 1) = A_g(p, 1)$ for all primes $p$.
Hence by the strong multiplicity one theorem (e.g.\ \cite[Theorem~12.6.1]{goldfeld2006}), we have $f = g$.

\section{Preliminaries}\label{sec:prelim}

Let $z = x + i y$ and let
\[
	E(z, s) = \frac{1}{2} \sum_{\substack{c, d \in \Z \\ (c, d) = 1}} \frac{y^s}{\abs{c z + d}^{2 s}}
\]
be the Eisenstein series with Fourier expansion
\[
	E(z, s) = y^s + \phi(s) y^{1 - s} + \sum_{n \neq 0} \phi(n, s) W_s(n z).
\]
Here $W_s$ is the Whittaker function
\[
	W_s(z) = 2 \abs{y}^{\frac{1}{2}} K_{s - \frac{1}{2}}(2 \pi \abs{y}) e(x),
\]
where $e(x) = e^{2 \pi i x}$, $K_s$ is the $K$-Bessel function,
\[
	\phi(s) = \sqrt{\pi} \frac{\Gamma\Bigl(s - \dfrac{1}{2}\Bigr)}{\Gamma(s)} \frac{\zeta(2 s - 1)}{\zeta(2 s)},
\]
\[
	\phi(n, s) = \pi^s \Gamma(s)^{-1} \zeta(2 s)^{-1} \abs{n}^{-\frac{1}{2}} \eta(n, s),
\]
and
\[
	\eta(n, s) = \sum_{a d = \abs{n}} \Bigl( \frac{a}{d} \Bigr)^{s - \frac{1}{2}}.
\]

For a Hecke--Maass form $f$ and the Eisenstein series $E = E(z, \frac{1}{2} + i t)$ we associate the Rankin--Selberg $L$-function
\[
		L(s, f \times E) = \sum_{m \geq 1} \sum_{n \geq 1} \frac{\bar \eta\Bigl(n, \dfrac{1}{2} + i t\Bigr) A(n, m)}{(m^2 n)^s}
\]
for $\Re(s) > 1$.
One derives
\[
	L\Bigl( \frac{1}{2}, f \times E\Bigr) = \abs[\Big]{L\Bigl( \frac{1}{2} - i t, f\Bigr) }^2,
\]
see e.g.\ \cite[Theorem~12.3.6]{goldfeld2006}.

For $u_j$ in the basis of Hecke--Maass forms for $\SL(2, \Z)$, we have the Fourier expansion
\[
	u_j(z) = \sum_{n \neq 0} \rho_j(n) W_{s_j}(n z).
\]
The normalization we use to get $\lambda_j(n)$ is
\[
	\rho_j(\pm n) = \rho_j(\pm 1) \lambda_j(n) n^{-\frac{1}{2}}.
\]

From Rankin--Selberg theory one gets
\[
	\mathop{\sum\sum}_{m^2 n \leq N} \abs{A(m, n)}^2 \ll N,
\]
which using Cauchy's inequality gives us
\[
	\sum_{n \leq N} \abs{A(m, n)} \ll N \abs{m}.
\]

\subsection{Kuznetsov trace formula}

For $m, n \geq 1$ and any even test function $h(t)$ holomorphic in the strip $\abs{\Im t} \leq \frac{1}{2} + \varepsilon$ and $h(t) \ll (\abs{t} + 1)^{-2 - \varepsilon}$ in the same strip, we have the Kuznetsov trace formula for even forms (see e.g., \cite[Section~3]{conreyiwaniec2000})
\begin{align*}
  &\sideset{}{'} \sum_j h(t_j) \omega_j \lambda_j(m) \lambda_j(n) + \frac{1}{4 \pi} \int_{-\infty}^\infty h(t) \omega(t) \conjugate{\eta}\Bigl ( m, \frac{1}{2} + i t \Bigr ) \eta\Bigl ( n, \frac{1}{2} + i t \Bigr ) \, d t \\
  &\quad= \frac{1}{2} \delta(m, n) H + \sum_{c > 0} \frac{1}{2 c} \left ( S(m, n; c) H^+\Bigl ( \frac{4 \pi \sqrt{m n}}{c} \Bigr ) + S(-m, n; c) H^- \Bigl ( \frac{4 \pi \sqrt{m n}}{c} \Bigr ) \right ).
\end{align*}
Here $'$ restricts the sum to even Maass forms, $\delta(m, n)$ is the Kronecker symbol,
\[
	\omega_j = \frac{4 \pi \abs{\rho_j(1)}^2}{\cosh(\pi t_j)},
\]
\[
	\omega(t) = \frac{4 \pi \abs[\Big]{\phi\Bigl(1, \dfrac{1}{2} + i t\Bigr)}^2}{\cosh(\pi t)}.
\]
\[
	H = \frac{2}{\pi} \int_0^\infty h(t) \tanh(\pi t) t \, d t,
\]
\[
	H^+(x) = 2 i \int_{-\infty}^\infty J_{2 i t}(x) \frac{h(t) t}{\cosh(\pi t)} \, d t,
\]
\[
	H^-(x) = \frac{4}{\pi} \int_{-\infty}^\infty K_{2 i t}(x) \sinh(\pi t) h(t) t \, d t,
\]
and
\[
	S(a, b; c) = \sum_{d \bar d \equiv 1 \pmod*{c}} e\Bigl( \frac{d a + \bar d b}{c} \Bigr)
\]
is the classical Kloosterman sum, along with $J_\nu$ and $K_\nu$ being the standard $J$-Bessel function and $K$-Bessel function respectively.

\subsection{Voronoi formula for $\GL(3)$}

Let
\[
	\alpha = -\nu_1 - 2 \nu_2 + 1, \quad \beta = -\nu_1 + \nu_2, \quad \gamma = 2 \nu_1 + \nu_2 - 1
\]
denote the Langlands parameters of $f$.

Per Miller and Schmid \cite[Theorem~1.18]{miller2006}, the Voronoi formula for $\GL(3)$ says that for $\psi$ a smooth compactly supported function on $(0, \infty)$, $c, d, \bar d \in \Z$ with $c \neq 0$, $(c, d) = 1$, and $d \bar d \equiv 1 \pmod{c}$, we have
\begin{align*}
	\sum_{n \geq 1} &A(n, m) e\Bigl( \frac{n d}{c} \Bigr) \psi(n) \\
	&= c \sum_{\pm} \sum_{n_1 \mid c m} \sum_{n_2 \geq 1} \frac{A(n_1, n_2)}{n_1 n_2} S(m \bar d, \pm n_2; m c n_1^{-1}) \Psi^{\pm} \Bigl( \frac{n_2 n_1^2}{c^3 m} \Bigr)\numberthis{\label{voronoiformula}}
\end{align*}
where
\[
	\Psi^\pm(x) = \int_{(\sigma)} y^{-s} \gamma^{\pm}(s) \tilde g(-s) \frac{d s}{2 \pi i}
\]
and $\sigma > \max\Set{-1 - \Re(\alpha), -1 - \Re(\beta), -1 - \Re(\gamma)}$.
Here
\[
	\tilde g(s) = \int_0^\infty x^{s - 1} g(x) \, d x
\]
is the Mellin transform of $g$ and
\[
	\gamma^{\pm}(s) = \gamma_0(s) \mp \gamma_1(s)
\]
with
\[
	\gamma_\ell(s) \coloneqq \frac{\pi^{-3 s - \frac{3}{2}}}{2} \frac{\Gamma\Bigl(\dfrac{1 + s + \alpha + \ell}{2}\Bigr) \Gamma\Bigl(\dfrac{1 + s + \beta + \ell}{2}\Bigr) \Gamma\Bigl(\dfrac{1 + s + \gamma + \ell}{2}\Bigr)}{\Gamma\Bigl(\dfrac{-s - \alpha + \ell}{2}\Bigr) \Gamma\Bigl(\dfrac{-s - \beta + \ell}{2}\Bigr) \Gamma\Bigl(\dfrac{-s - \gamma + \ell}{2}\Bigr)}
\]
for $\ell = 0, 1$.

From \cite[Lemma~2.1]{li2011}, if $\psi$ in addition to being smooth is compactly supported on $[X, 2 X]$, then for any fixed integer $K \geq 1$ and $x X \gg 1$, we have
\begin{equation}\label{voronoiformulaasymptotics}
	\Psi^\pm(x) = x \int_0^\infty \psi(y) \sum_{j = 1}^K \frac{c_j^\pm e(3 x^\frac{1}{3} y^{\frac{1}{3}}) + d_j^\pm e(-3 x^{\frac{1}{3}} y^{\frac{1}{3}})}{(x y)^{\frac{j}{3}}} \, d y + O\bigl( (x X)^{\frac{-K + 2}{3}} \bigr)
\end{equation}
where $c_j^\pm$ and $d_j^\pm$ are absolute constants depending on $\alpha$, $\beta$, and $\gamma$.
In practice we will only work with the $j = 1$ term coming from this, since all other terms are lower order and hence smaller and give even better estimates.

\subsection{Approximate functional equations}

Now $L(s, f \times u_j)$ has an approximate functional equation (see \cite[Theorem~5.3]{iwaniec2014}), namely
\begin{align*}
  L\Bigl( \frac{1}{2}, f \times u_j\Bigr ) &= \sum_{m \geq 1} \sum_{n \geq 1} \frac{\lambda_j(n) A(n, m)}{(m^2 n)^{\frac{1}{2}}} V_-(m^2 n, t_j) + {} \\
	&\qquad\qquad {} + \sum_{m \geq 1} \sum_{n \geq 1} \frac{\lambda_j(n) A(m, n)}{(m^2 n)^{\frac{1}{2}}} V_+(m^2 n, t_j)
\end{align*}
where for
\[
  F(u) = \Bigl ( \cos \frac{\pi u}{A} \Bigr)^{-3 A},
\]
with $A$ a positive integer, and for $\abs{\Im t} \leq 1000$, we have
\[
  V_\mp(y, t) = \frac{1}{2 \pi i} \int_{(1000)} y^{- u} F(u) \frac{\gamma_\mp\Bigl(\dfrac{1}{2} + u, t\Bigr)}{\gamma_-\Bigl(\dfrac{1}{2}, t\Bigr)} \frac{d u}{u}
\]
and
\begin{gather*}
	\gamma_\mp(s, t) = \pi^{-3 s} \Gamma\Bigl ( \frac{s - i t \mp \alpha}{2} \Bigr ) \Gamma\Bigl ( \frac{s - i t \mp \beta}{2} \Bigr ) \Gamma\Bigl ( \frac{s - i t \mp \gamma}{2} \Bigr ) \times{} \\
	\qquad{} \times \Gamma\Bigl ( \frac{s + i t \mp \alpha}{2} \Bigr ) \Gamma\Bigl ( \frac{s + i t \mp \beta}{2} \Bigr ) \Gamma\Bigl ( \frac{s + i t \mp \gamma}{2} \Bigr ). \numberthis\label{gammafactors}
\end{gather*}
Following \cite[Proposition~5.4]{iwaniec2014}, essentially by Stirling's formula, the growth of $V_\pm(m^2n, t_j)$ limit the sums in $L(\frac{1}{2}, f \times u_j)$ to $m^2 n \ll \abs{t_j}^{3 + \varepsilon}$.

In order to make use of the Kuznetsov trace formula we also need an approximate functional equation for $L(\frac{1}{2}, f \times E) = \abs{L(\frac{1}{2} - i t, f)}^2$.
As can be verified by using the functional equation of $L(s, f)$ itself, this Rankin--Selberg $L$-function satisfies the same functional equation as $L(s, f \times u_j)$ and consequently
\begin{align*}
  \abs[\Big]{L\Bigl ( \frac{1}{2} - i t, f \Bigr)}^2 &= \sum_{m \geq 1} \sum_{n \geq 1} \frac{\eta\Bigl(n, \dfrac{1}{2} + i t\Bigr) A(n, m)}{(m^2 n)^{\frac{1}{2}}} V_-(m^2 n, t) + {} \\
	&\qquad\qquad {} + \sum_{m \geq 1} \sum_{n \geq 1} \frac{\eta\Bigl(n, \dfrac{1}{2} + i t\Bigr) A(m, n)}{(m^2 n)^{\frac{1}{2}}} V_+(m^2 n, t).
\end{align*}

\section{Proof of \autoref{thm:maintheorem}: The setup}\label{sec:setup}

Substituting these approximate functional equations into the spectrally normalized moment we wish to compute we get
\begin{align*}
  &\sum_{m \geq 1} \sum_{n \geq 1} \frac{A(n, m)}{(m^2 n)^{\frac{1}{2}}} \biggl ( \sideset{}{'} \sum_j k(t_j) \omega_j \lambda_j(p) \lambda_j(n) V_-(m^2 n, t_j) + {} \\
  &\qquad\qquad{}+ \frac{1}{4 \pi} \int_{-\infty}^\infty k(t) \omega(t) \conjugate{\eta}\Bigl ( p, \frac{1}{2} + i t \Bigr ) \eta\Bigl ( n, \frac{1}{2} + i t \Bigr ) V_-(m^2 n, t) \, d t \biggr ) + {} \\
	&{} + \sum_{m \geq 1} \sum_{n \geq 1} \frac{A(m, n)}{(m^2 n)^{\frac{1}{2}}} \biggl ( \sideset{}{'} \sum_j k(t_j) \omega_j \lambda_j(p) \lambda_j(n) V_+(m^2 n, t_j) + {} \\
  &\qquad\qquad{}+ \frac{1}{4 \pi} \int_{-\infty}^\infty k(t) \omega(t) \conjugate{\eta}\Bigl ( p, \frac{1}{2} + i t \Bigr ) \eta\Bigl ( n, \frac{1}{2} + i t \Bigr ) V_+(m^2 n, t) \, d t \biggr ).
\end{align*}
Hence computing asymptotics for this spectrally normalized moment proves \autoref{thm:maintheorem}.

Because $V_-(y, t)$ and $V_+(y, t)$ are essentially the same, we will suppress the subscript in the following discussion and deal with the first sum above, keeping in mind that the dual sum behaves exactly the same, only its Fourier coefficient has the arguments reversed.

By calling $h(t) = k(t) V(m^2 n, t)$ we can apply our Kuznetsov trace formula to the first sum, resulting in
\begin{gather*}
  \frac{1}{2} \sum_{m \geq 1} \sum_{n \geq 1} \frac{A(n, m)}{(m^2 n)^{\frac{1}{2}}} \Biggl ( \delta(n, p) H_{m, n} + \sum_{c > 0} \frac{1}{c} \biggl (  S(n, p; c) H_{m, n}^+\Bigl ( \frac{4 \pi \sqrt{n p}}{c} \Bigr ) + {} \\
  {}+S(-n, p; c) H_{m, n}^- \Bigl ( \frac{4 \pi \sqrt{n p}}{c} \Bigr ) \biggr ) \Biggr ),
\end{gather*}
where this time in particular
\[
  H_{m, n} = \frac{2}{\pi} \int_0^\infty k(t) V(m^2 n, t) \tanh(\pi t) t \, d t,
\]
\[
	H_{m, n}^+(x) = 2 i \int_{-\infty}^\infty J_{2 i t}(x) \frac{k(t) V(m^2 n, t) t}{\cosh(\pi t)} \, d t,
\]
and
\[
	H_{m, n}^-(x) = \frac{4}{\pi} \int_{-\infty}^\infty K_{2 i t}(x) \sinh(\pi t) k(t) V(m^2 n, t) t \, d t.
\]
We split the resulting sum into
\[
	\mathcal{D} + \mathcal{O}^+ + \mathcal{O}^-
\]
where
\[
	\mathcal{D} = \frac{1}{2} \sum_{m \geq 1} \sum_{n \geq 1} \frac{A(n, m)}{(m^2 n)^{\frac{1}{2}}} \delta(n, p) H_{m, n},
\]
\[
	\mathcal{O^+} = \frac{1}{2} \sum_{m \geq 1} \sum_{n \geq 1} \frac{A(n, m)}{(m^2 n)^{\frac{1}{2}}} \sum_{c > 0} \frac{S(n, p; c)}{c} H_{m, n}^+ \Bigl (\frac{4 \pi \sqrt{n p}}{c} \Bigr ),
\]
and
\[
	\mathcal{O}^- = \frac{1}{2} \sum_{m \geq 1} \sum_{n \geq 1} \frac{A(n, m)}{(m^2 n)^{\frac{1}{2}}} \sum_{c > 0} \frac{S(-n, p; c)}{c} H_{m, n}^- \Bigl ( \frac{4 \pi \sqrt{n p}}{c} \Bigr ).
\]
We will estimate each of these three terms in turn.

\section{Proof of \autoref{thm:maintheorem}: The diagonal terms}\label{sec:diagonal}

We prove the following:

\begin{lemma}
	\begin{align*}
		\mathcal{D} &= \frac{1}{2} \sum_{m \geq 1} \sum_{n \geq 1} \frac{A(n, m)}{(m^2 n)^{\frac{1}{2}}} \delta(n, p) H_{m, n} \\
		&= \frac{L(1, \tilde f) \bigl(A(p, 1) p - 1\bigr)}{p^{\frac{3}{2}} \pi} \int_0^\infty k(t) \tanh(\pi t) t \, d t + O(M T^{\frac{1}{7} + \varepsilon} p^{\varepsilon}).
	\end{align*}
\end{lemma}
The dual sum is identical, only with $L(1, f)$ instead of $L(1, \tilde f)$ and $A(1, p)$ instead of $A(p, 1)$.

Since $\delta(n, p) = 1$ if $n = p$ and $0$ otherwise, we have
\[
  \mathcal{D} = \frac{1}{2} \sum_{m \geq 1} \frac{A(p, m)}{m p^{\frac{1}{2}}} H_{m, p}.
\]
We can write this sum as
\begin{align*}
  \mathcal{D} &= \frac{1}{2} \sum_{(m, p) = 1} \frac{A(p, m)}{m p^{\frac{1}{2}}} H_{m, p} + \frac{1}{2} \sum_{(m, p) > 1} \frac{A(p, m)}{m p^{\frac{1}{2}}} H_{m, p} \\
  &= \frac{1}{2} \sum_{(m, p) = 1} \frac{A(p, m)}{m p^{\frac{1}{2}}} H_{m, p} + \frac{1}{2} \sum_{m \geq 1} \frac{A(p, m p)}{(m p) p^{\frac{1}{2}}} H_{m p, p}.
\end{align*}
For $(m, p) = 1$, the multiplicativity of the Fourier coefficients gives us $A(p, m) = A(p, 1) A(1, m)$, and in the second sum by the Hecke relations (see e.g., \cite[Theorem~6.4.11]{goldfeld2006})
\begin{align*}
  A(p, m p) &= \sum_{d \mid (p, p m)} \mu(d) A\Big ( \frac{p}{d}, 1 \Big ) A \Big ( 1, \frac{p m}{d} \Big ) \\
  &= A(p, 1) A(1, p m) - A(1, 1) A(1, m) \\
  &= A(p, 1) A(1, p m) - A(1, m)
\end{align*}
since $A(1, 1) = 1$ by normalization.
Putting these expressions back into $\mathcal{D}$ and rearranging we get
\begin{equation}\label{diagonalterms}
  \mathcal{D} = \frac{A(p, 1)}{2 p^{\frac{1}{2}}} \sum_{m \geq 1} \frac{A(1, m)}{m} H_{m, p} - \frac{1}{2 p^{\frac{3}{2}}} \sum_{m \geq 1} \frac{A(1, m)}{m} H_{m p, p}.
\end{equation}

It remains to estimate the two sums.
Recalling now
\[
	H_{m, p} = \frac{2}{\pi} \int_0^\infty k(t) V(m^2 p, t) \tanh(\pi t) t \, d t
\]
and
\[
	V(m^2 p, t) = \frac{1}{2 \pi i} \int_{(1000)} (m^2 p)^{- u} F(u) \frac{\gamma\Bigl(\dfrac{1}{2} + u, t\Bigr)}{\gamma\Bigl(\dfrac{1}{2}, t\Bigr)} \frac{d u}{u},
\]
we can bring the sum all the way into $V(m^2 p, t)$, getting
\begin{align*}
	&\sum_{m \geq 1} \frac{A(1, m)}{m} H_{m, p} \\
	&\quad= \frac{2}{\pi} \int_0^\infty k(t) \tanh(\pi t) t \biggl ( \frac{1}{2 \pi i} \int_{(1000)} p^{-u} \Bigl ( \sum_{m \geq 1} \frac{A(1, m)}{m^{1 + 2 u}}\Bigr ) F(u) \frac{\gamma\Bigl(\dfrac{1}{2} + u, t\Bigr)}{\gamma\Bigl(\dfrac{1}{2}, t\Bigr)} \frac{d u}{u} \biggr ) \, d t.
\end{align*}
The sum on the inside is just $L(1 + 2 u, \tilde f)$ (keep in mind this means the dual sum will have $L(1 + 2 u, f)$ here instead).

By the Kim--Sarnak bound on the generalized Ramanujan conjecture for $\GL(n)$, see \cite[Appendix~2]{Kim2003},
\[
	\abs{\Re(\alpha)}, \abs{\Re(\beta)}, \abs{\Re(\gamma)} \leq \frac{1}{2} - \frac{1}{7}.
\]
Hence we can shift the line of integration to $d = -\frac{1}{7} + \varepsilon$ without hitting any of the poles of the gamma factors.
Doing this we pass the simple pole at $u = 0$, picking up the residue $L(1, \tilde f)$.
Therefore the inner integral becomes
\begin{equation}\label{eq:innerintegral-mainterm}
	L(1, \tilde f) + \frac{1}{2 \pi i} \int_{(d)} p^{-u} L(1 + 2 u, \tilde f) F(u) \frac{\gamma\Bigl(\dfrac{1}{2} + u, t\Bigr)}{\gamma\Bigl(\dfrac{1}{2}, t\Bigr)} \frac{d u}{u}.
\end{equation}

To estimate the inner integral we need to control the gamma factors in terms of $t$.
By Stirling's formula (see also \cite[Proposition~5.4]{iwaniec2014}) we have
\[
	\frac{\Gamma(s + u)}{\Gamma(s)} \ll \abs{s + 1}^{\Re(u)} \exp\Bigl( \frac{\pi}{2} \abs{u} \Bigr)
\]
for $\Re(s) > 0$ and $\Re(u) > -\Re(s)$, so in $t$ each gamma factor contributes $t^{-\frac{1}{7} + \varepsilon}$.

The resulting exponential in $\abs{u}$ is controlled by $F(u)$ in the integral, and consequently what we get out of the inner integral \eqref{eq:innerintegral-mainterm} is $L(1, \tilde f) + O(p^{\frac{1}{7} - \varepsilon} t^{-\frac{6}{7} + \varepsilon})$, since we have six such gamma factors, and from this we get
\begin{align*}
	\sum_{m \geq 1} \frac{A(1, m)}{m} H_{m, p} &= L(1, \tilde f) \frac{2}{\pi} \int_0^\infty k(t) \tanh(\pi t) t \, d t + {} \\
	&\qquad{} + O \Bigl ( p^{\frac{1}{7} - \varepsilon} \int_0^\infty k(t) \tanh(\pi t) t^{\frac{1}{7} + \varepsilon} \, d t \Bigr ).\numberthis\label{mainterm}
\end{align*}
The second term in \eqref{diagonalterms} is treated similarly, the only difference being a factor of $p^{-3 u}$ in the integral instead of $p^{-u}$, hence
\begin{align*}
	\sum_{m \geq 1} \frac{A(1, m)}{m} H_{m p, p} &= L(1, \tilde f) \frac{2}{\pi} \int_0^\infty k(t) \tanh(\pi t) t \, d t + {} \\
	&\qquad{} + O \Bigl ( p^{\frac{3}{7} - \varepsilon} \int_0^\infty k(t) \tanh(\pi t) t^{\frac{1}{7} + \varepsilon} \, d t \Bigr ).
\end{align*}

Putting both of these back into \eqref{diagonalterms} we get
\begin{align*}
	\mathcal{D} &= \frac{L(1, \tilde f) \bigl(A(p, 1) p - 1\bigr)}{p^{\frac{3}{2}} \pi} \int_0^{\infty} k(t) \tanh(\pi t) t \, d t + O(M T^{\frac{1}{7} + \varepsilon} p^{\varepsilon})
\end{align*}
since $t$ is about $T$ and the integral is of length about $M$.
Here we have again used the Kim--Sarnak bound, this time in the form $\abs{A(p, 1)} \ll p^{\frac{5}{14} + \varepsilon}$ (see \cite[Appendix~2]{Kim2003}) in the error term.

Working through the exact same calculations for the dual sum we consequently get
\begin{align*}
	&\frac{A(1, p)}{2 p^{\frac{1}{2}}} \sum_{m \geq 1} \frac{A(m, 1)}{m} H_{m, p} - \frac{1}{2 p^{\frac{3}{2}}} \sum_{m \geq 1} \frac{A(m, 1)}{m} H_{m p, p} \\
	&\qquad= \frac{L(1, f) \bigl(A(1, p) p - 1\bigr)}{p^{\frac{3}{2}} \pi} \int_0^{\infty} k(t) \tanh(\pi t) t \, d t + O(M T^{\frac{1}{7} + \varepsilon} p^{\varepsilon})
\end{align*}
which means the two diagonal sums put together contribute
\begin{align*}
	&\frac{L(1, \tilde f) \bigl(A(p, 1) p - 1\bigr) + L(1, f) \bigl(A(1, p) p - 1\bigr)}{p^{\frac{3}{2}} \pi} \int_0^{\infty} k(t) \tanh(\pi t) t \, d t + {} \\
	&\qquad{} + O(M T^{\frac{1}{7} + \varepsilon} p^{\varepsilon})
\end{align*}
to the moment in \autoref{thm:maintheorem}.

\section{Proof of \autoref{thm:maintheorem}: The $J$-Bessel function terms}\label{sec:Jbessel}

For $\mathcal{O}^+$, using partition of unity it suffices to consider
\[
	\mathcal{R^+} = \frac{1}{2} \sum_{m \geq 1} \sum_{n \geq 1} \frac{A(n, m)}{(m^2 n)^{\frac{1}{2}}} g\Bigl( \frac{m^2 n}{N} \Bigr) \sum_{c > 0} \frac{S(n, p; c)}{c} H_{m, n}^+ \Bigl (\frac{4 \pi \sqrt{n p}}{c} \Bigr )
\]
where $g$ is a smooth function of compact support on $\interval{1}{2}$ and $N$ is at most $T^{3 + \varepsilon}$.

We prove the following:
\begin{lemma}
	\begin{align*}
		\mathcal{R^+} = O(M^{-3} T^{\frac{5}{2} + \varepsilon} p^{1 + \varepsilon}).
	\end{align*}
\end{lemma}

Since the off-diagonal terms will contribute only error terms we will have no need to keep track of the difference between this sum and its dual version.
In what follows we will again suppress the subscripts on $V_\pm(m^2 n, t)$, but the dual sum contributes precisely the same error.

Following \cite{li2011}, the strategy for estimating $\mathcal{R}^+$ is to split the $c$-sum into three parts,
\[
	\mathcal{R}^+ = \mathcal{R}^+_1 + \mathcal{R}^+_2 + \mathcal{R}^+_3,
\]
with
\[
	\mathcal{R}^+_1 = \frac{1}{2} \sum_{m \geq 1} \sum_{n \geq 1} \frac{A(n, m)}{(m^2 n)^{\frac{1}{2}}} g\Bigl( \frac{m^2 n}{N} \Bigr) \sum_{c \geq \frac{C_1}{m}} \frac{S(n, p; c)}{c} H_{m, n}^+ \Bigl( \frac{4 \pi \sqrt{n p}}{c} \Bigr),
\]
\[
	\mathcal{R}^+_2 = \frac{1}{2} \sum_{m \geq 1} \sum_{n \geq 1} \frac{A(n, m)}{(m^2 n)^{\frac{1}{2}}} g\Bigl( \frac{m^2 n}{N} \Bigr) \sum_{\frac{C_2}{m} \leq c \leq \frac{C_1}{m}} \frac{S(n, p; c)}{c} H_{m, n}^+ \Bigl( \frac{4 \pi \sqrt{n p}}{c} \Bigr),
\]
and
\[
	\mathcal{R}^+_3 = \frac{1}{2} \sum_{m \geq 1} \sum_{n \geq 1} \frac{A(n, m)}{(m^2 n)^{\frac{1}{2}}} g\Bigl( \frac{m^2 n}{N} \Bigr) \sum_{c \leq \frac{C_2}{m}} \frac{S(n, p; c)}{c} H_{m, n}^+ \Bigl( \frac{4 \pi \sqrt{n p}}{c} \Bigr).
\]
The idea is to tune the cut-offs $C_1$ and $C_2$ in such a way that both $\mathcal{R}^+_1$ and $\mathcal{R}^+_2$ are small, meaning negative powers of $T$, and $\mathcal{R}^+_3$ we will handle a bit more delicately with stationary phase analysis.

Starting with the tail $\mathcal{R}^+_1$, we move the line of integration in
\[
	H_{m, n}^+ (x) = 2 i \int_{-\infty}^\infty J_{2 i t}(x) \frac{k(t) V(m^2 n, t) t}{\cosh(\pi t)} \, d t
\]
to $\Im t = -\frac{1}{2}$.
This gives us
\[
	H_{m, n}^+ (x) = 2 i \int_{-\infty}^\infty J_{2 i y + 1}(x) \frac{k\Bigl(- \dfrac{1}{2} i + y\Bigr) V\Bigl(m^2 n, - \dfrac{1}{2} i + y\Bigr) \Bigl(- \dfrac{1}{2} i + y\Bigr)}{\cosh\Bigl(\pi \Bigl(- \dfrac{1}{2} i + y\Bigr)\Bigr)} \, d y
\]
since this shift passes through no zeros of the denominator $\cosh(\pi t)$.

The $J$-Bessel function has the following integral representation (see for instance \cite[8.411~4]{gradshtein2015})
\[
	J_{\nu}(z) = 2 \frac{\Bigl(\dfrac{z}{2}\Bigr)^\nu}{\Gamma\Bigl(\nu + \dfrac{1}{2}\Bigr) \Gamma\Bigl(\dfrac{1}{2}\Bigr)} \int_0^{\pi/2} \sin(\theta)^{2 \nu} \cos(z \cos \theta) \, d \theta,
\]
which tells us that for $\Re \nu > -\frac{1}{2}$,
\[
	J_\nu(z) \ll \Bigl( \frac{z}{\abs{\Im \nu} + 1} \Bigr)^{\Re \nu} e^{\frac{\pi}{2} \abs{\Im \nu} }.
\]
Using Stirling's formula to estimate the gamma factors in $V(m^2 n, -d i + y)$, we get
\begin{equation}\label{eq:StirlingV}
	V(m^2 n, -d i + y) \ll \Bigl ( \frac{\abs{y + 1}^3}{m^2 n} \Bigr)^d.
\end{equation}

Putting this together and back into $H_{m, n}^+(x)$, keeping in mind that this integral is of length about $M$ because of the exponential decay of $k(-\frac{1}{2} i + y)$, and that $y$ is about $T$, this gives us
\[
	H_{m, n}^+(x) \ll x T^{-1} (m^2 n)^{-\frac{1}{2}} T^{\frac{3}{2}} T^{1 + \varepsilon} M = x T^{\frac{3}{2} + \varepsilon} (m^2 n)^{-\frac{1}{2}} M.
\]

Plugging this back into $\mathcal{R}^+_1$, summing trivially over $c$ using Weil's bound
\[
	S(a, b; c)  \ll_\varepsilon c^{\frac{1}{2} + \varepsilon} (a, b, c)^{\frac{1}{2}}
\]
for the Kloosterman sum, this gives us
\begin{align*}
	\sum_{c \geq \frac{C_1}{m}} \frac{S(n, p; c)}{c} H_{m, n}^+ \Bigl( \frac{4 \pi \sqrt{n p}}{c} \Bigr) &\ll \sum_{c \geq \frac{C_1}{m}} \frac{c^{\frac{1}{2} + \varepsilon}}{c} \Bigl( \frac{4 \pi \sqrt{n p}}{c} \Bigr) T^{\frac{3}{2} + \varepsilon} (m^2 n)^{-\frac{1}{2}} M \\
	&= (4 \pi) (n p)^{\frac{1}{2}} T^{\frac{3}{2} + \varepsilon} (m^2 n)^{-\frac{1}{2}} M \sum_{c \geq \frac{C_1}{m}} c^{\frac{1}{2} - 1 - 1 + \varepsilon} \\
	&= (4 \pi) T^{\frac{3}{2} + \varepsilon} m^{-1} p^{\frac{1}{2}} M \sum_{c \geq \frac{C_1}{m}} c^{-\frac{3}{2} + \varepsilon}.
\end{align*}
Hence the $c$-sum converges.
In particular, the sum over $c$ is about $( \frac{C_1}{m} )^{-\frac{1}{2} + \varepsilon}$.
By picking $C_1$ to be an appropriately large power of $T$, say $T^{100}$, we can ensure that the power of $T$ is $\frac{3}{2} + \varepsilon + 100 ( - \frac{1}{2} + \varepsilon ) < 0$.
Hence any power of $T$ larger than $3$ will do.
In other words, $\mathcal{R}^+_1$ contributes a negative power of $T$, and as such is admissible.

Moving on to
\[
	\mathcal{R}^+_2 = \frac{1}{2} \sum_{m \geq 1} \sum_{n \geq 1} \frac{A(n, m)}{(m^2 n)^{\frac{1}{2}}} g\Bigl( \frac{m^2 n}{N} \Bigr) \sum_{\frac{C_2}{m} \leq c \leq \frac{C_1}{m}} \frac{S(n, p; c)}{c} H_{m, n}^+ \Bigl( \frac{4 \pi \sqrt{n p}}{c} \Bigr)
\]
our goal is to choose $C_2$ in such a way that this, too, is a negative power of $T$, keeping in mind that $C_1 = T^{100}$.
To do this we will use another integral representation of the $J$-Bessel function, namely
\[
	\frac{J_{2 i t}(x) - J_{-2 i t}(x)}{\cosh(\pi t)} = - \frac{2 i}{\pi} \tanh(\pi t) \int_{-\infty}^\infty \cos(x \cosh \zeta) e\Bigl( \frac{t \zeta}{\pi} \Bigr) \, d \zeta,
\]
derived from \cite[8.411~11]{gradshtein2015}.
In other words we want to decompose $H_{m, n}(x)$ as follows:
\begin{align*}
	H_{m, n}^+(x) &= 2 i \int_{-\infty}^\infty J_{2 i t}(x) \frac{k(t) V(m^2 n, t) t}{\cosh(\pi t)} \, d t \\
	&= 2 i \int_{0}^\infty J_{2 i t}(x) \frac{k(t) V(m^2 n, t) t}{\cosh(\pi t)} \, d t + 2 i \int_{-\infty}^0 J_{2 i t}(x) \frac{k(t) V(m^2 n, t) t}{\cosh(\pi t)} \, d t \\
	&= 2 i \int_{0}^\infty J_{2 i t}(x) \frac{k(t) V(m^2 n, t) t}{\cosh(\pi t)} \, d t - 2 i \int_{0}^\infty J_{-2 i t}(x) \frac{k(-t) V(m^2 n, -t) t}{\cosh(-\pi t)} \, d t.
\end{align*}
But $k(t)$, $V(m^2n, t)$, and $\cosh(\pi t)$ are all even in $t$, so this becomes
\[
	H_{m, n}^+(x) = 2 i \int_0^\infty \frac{J_{2 i t}(x) - J_{- 2 i t}(x)}{\cosh(\pi t)} k(t) V(m^2 n, t) t \, d t,
\]
and so we can apply the above integral representation to get
\[
	H_{m, n}^+(x) = \frac{4}{\pi} \int_{t = 0}^\infty \tanh(\pi t) \biggl ( \int_{\zeta = -\infty}^\infty \cos(x \cosh \zeta) e\Bigl( \frac{t \zeta}{\pi} \Bigr) \, d \zeta \biggr ) k(t) V(m^2n, t) t \, d t.
\]
We can get rid of $e^{-\frac{(t + T)^2}{M^2}}$ from $k(t) = e^{-\frac{(t - T)^2}{M^2}} + e^{-\frac{(t + T)^2}{M^2}}$ with negligible error, and similarly $\tanh(\pi t)$ is inconsequential, so we are left with studying
\[
	H_{m, n}^+(x) = \frac{4}{\pi} \int_{t = 0}^\infty \int_{\zeta = -\infty}^\infty \cos(x \cosh \zeta) e\Bigl( \frac{t \zeta}{\pi} \Bigr) e^{- \frac{(t - T)^2}{M^2}} V(m^2 n, t) t \, d \zeta \, d t + O(T^{-A})
\]
where $A$ is an arbitrarily large constant.
We make the change of variables $u = \frac{t - T}{M}$, under which our integral becomes
\begin{align*}
	H_{m, n}^+(x) &= \frac{4 M}{\pi} \int_{u = -\frac{T}{M}}^\infty \int_{\zeta = -\infty}^\infty \cos(x \cosh \zeta) e\Bigl ( \frac{(M u + T) \zeta}{\pi} \Bigr) (M u + T) e^{-u^2} \times{} \\
	&\qquad\qquad\qquad {}\times V(m^2 n, M u + T) \, d \zeta \, d u + O(T^{-A}).
\end{align*}
We only study the term arising from the $T$ part of $M u + T$; the term from $M u$ can be handled similarly.

Hence we want to understand
\[
	H_{m, n}^{+, 1}(x) = \frac{4 M T}{\pi} \int_{u = -\infty}^\infty \int_{\zeta = -\infty}^\infty e^{-u^2} V(m^2 n, M u + T) \cos(x \cosh \zeta) e\Bigl( \frac{u M \zeta}{\pi} \Bigr) e\Bigl( \frac{T \zeta}{\pi} \Bigr) \, d u \, d \zeta,
\]
where we have extended the $u$-integral to $\interval[open]{-\infty}{\infty}$ with negligible error because of the exponential decay of $e^{-u^2}$.

By setting
\[
	k^*(u) = e^{-u^2} V(m^2 n, M u + T)
\]
and considering its Fourier transform
\[
	\widehat{k^*}(\zeta) = \int_{-\infty}^\infty k^*(u) e(-u \zeta) \, d u,
\]
we get in particular that
\[
	\widehat{k^*}\Bigl( -\frac{M \zeta}{\pi} \Bigr) = \int_{-\infty}^\infty k^*(u) e\Bigl( \frac{u M \zeta}{\pi} \Bigr) \, d u.
\]
This lets us rewrite $H_{m, n}^{+, 1}(x)$ as
\[
	H_{m, n}^{+, 1}(x) = \frac{4 M T}{\pi} \int_{-\infty}^\infty \widehat{k^*}\Bigl( -\frac{M \zeta}{\pi} \Bigr) \cos(x \cosh \zeta) e\Bigl( \frac{T \zeta}{\pi} \Bigr) \, d \zeta,
\]
and changing variables to $\xi = - \frac{M \zeta}{\pi}$, this leaves us studying
\begin{equation}\label{eq:Hmnplus}
	H_{m, n}^{+, 1}(x) = 4 T \int_{-\infty}^\infty \widehat{k^*}(\xi) \cos\Bigl ( x \cosh \frac{\xi \pi}{M} \Bigr) e \Bigl ( - \frac{T \xi}{M} \Bigr) \, d \xi.
\end{equation}
Rewriting cosine in terms of exponentials, the phase of this integral is
\[
	\phi(\xi) = - \frac{T \xi}{M} \pm \frac{x}{2 \pi} \cosh \frac{\xi \pi}{M},
\]
the plus or minus depending on which half of the exponential representation of cosine we are considering, but the analysis is the same either way:
\[
	\abs{\phi'(\xi)} \geq \frac{T}{M} - \frac{x}{2 M} \sinh \frac{\xi \pi}{M} \gg \frac{T}{M} - \frac{x}{2 M} \frac{\xi \pi}{M} = \frac{T}{M} - \frac{x \xi \pi}{2 M^2}
\]
Hence for there to be a stationary phase we must have $x$ about the size $T M$, and consequently so long as $\abs{x} \leq T^{1 - \varepsilon} M$, there is no stationary phase, and so by the First derivative test (\cite[Lemma~5.1.2]{huxley1996}) $H_{m, n}^{+, 1}(x)$ (and hence $H_{m, n}^+(x)$, being composed of $H_{m, n}^{+, 1}(x)$ and a similar term) is small in this region.

In particular, since $x = \frac{4 \pi \sqrt{n p}}{c}$, this means that $\mathcal{R}^+_2$ is negligible for $\frac{4 \pi \sqrt{n p}}{c} \leq T^{1 - \varepsilon} M$, or in other words
\[
	c \geq \frac{4 \pi \sqrt{n p}}{T^{1 - \varepsilon} M}.
\]
As promised we want to tune $C_2$ in such a way that $\mathcal{R}^+_2$, the $c$-sum of which ranges over $\frac{C_2}{m} \leq c \leq \frac{C_1}{m}$, is small.
Remembering that $1 \leq \frac{m^2 n}{N} \leq 2$, i.e., $n \gg N m^{-2}$, we infer from the above that
\[
	C_2 = \frac{\sqrt{N p}}{T^{1 - \varepsilon} M}
\]
does the job, since in this range the resulting exponential integral has no stationary phase, and recall that in $\mathcal{R}^+_2$, the $c$-sum is over $\frac{C_2}{m} \leq c \leq \frac{C_1}{m}$, hence convergent.

Finally for
\[
	\mathcal{R}^+_3 = \frac{1}{2} \sum_{m \geq 1} \sum_{n \geq 1} \frac{A(n, m)}{(m^2 n)^{\frac{1}{2}}} g\Bigl( \frac{m^2 n}{N} \Bigr) \sum_{c \leq \frac{C_2}{m}} \frac{S(n, p; c)}{c} H_{m, n}^+ \Bigl( \frac{4 \pi \sqrt{n p}}{c} \Bigr).
\]
we need to be substantially more careful, because this time, for $c \leq \frac{C_2}{m}$, the integral $H_{m, n}^{+, 1}$ in \eqref{eq:Hmnplus} has stationary phase for $\xi \asymp \pm \frac{2 M T}{\pi x}$.
We will focus on the stationary phase near $\xi_0 = -\frac{2 M T}{\pi x}$ from the negative part of $\pm$ in the phase $\phi(\xi)$ simply to keep the signs uniform; the positive one is treated similarly.
To this end we will call
\[
	\tilde H_{m, n}^{+, 1}(x) = 4 T \int_{-\infty}^\infty \widehat{k^*}(\xi) e\Bigl( -\frac{T \xi}{M} - \frac{x}{2 \pi} \cosh \frac{\xi \pi}{M} \Bigr) \, d \xi.
\]

Following \cite[Proposition~4.1]{li2011} (or \cite[Lemma~5.1]{lau2006} or \cite[Proposition~3.1]{sarnak2001}), we can extract asymptotics for $\tilde H_{m, n}^{+, 1}(x)$.
This is essentially ordinary stationary phase analysis, only being quite careful with the computation of the error.

By expanding the $\cosh$ in the phase $\phi(\xi)$ as a Taylor series we get
\begin{align*}
	\tilde H_{m, n}^{+, 1}(x) &= 4 T \int_{-\infty}^\infty \widehat{k^*}(\xi) e\Bigl(-\frac{T \xi}{M} - \frac{x}{2 \pi} - \frac{\pi x \xi^2}{4 M^2} - \frac{\pi^3 x \xi^4}{48 M^4} - \frac{\pi^5 x \xi^6}{1440 M^6} \Bigr) \, d \xi + {} \\
	& \qquad\qquad\qquad {} + O\Bigl( T \int_{-\infty}^\infty \abs{\widehat{k^*}(\xi)} \frac{\abs{\xi}^8 \abs{x}}{M^8} \, d \xi \Bigr).
\end{align*}
Note that the integral in the error is finite in $\xi$, so the error here is $O(\frac{T \abs{x}}{M^8})$.
By expanding the last exponential term $e(-\frac{\pi^5 x \xi^6}{1440 M^6})$ as a Taylor series of order $1$ we get
\begin{equation}\label{Hmndecomp}
	\tilde H_{m, n}^{+, 1}(x) = W_{m, n}^+(x) - \frac{2 \pi^6 x}{1440 M^6} W_{m, n}^-(x) + O\Bigl( \frac{T \abs{x}}{M^8} \Bigr)
\end{equation}
where
\[
	W_{m, n}^+(x) = 4 T e\Bigl(-\frac{x}{2 \pi}\Bigr) \int_{-\infty}^\infty k_0^*(\xi) e\Bigl( -\frac{T \xi}{M} - \frac{\pi x \xi^2}{4 M^2} \Bigr) \, d \xi
\]
with
\[
	k_0^*(\xi) = \widehat{k^*}(\xi) e\Bigl( -\frac{\pi^3 x \xi^4}{48 M^4} \Bigr)
\]
and
\[
	W_{m, n}^-(x) = 4 T e\Bigl( - \frac{x}{2 \pi} \Bigr) \int_{-\infty}^\infty \xi^6 k_0^*(\xi) e\Bigl( -\frac{T\xi}{M} - \frac{\pi x \xi^2}{4 M^2} \Bigr) \, d \xi.
\]
In what follows we will deal with $W_{m, n}^+(x)$ as $W_{m, n}^-(x)$ can be handled with precisely the same methods.
As detailed in the calculations above \cite[Proposition~4.1]{li2011}, by completing the square in the exponential, applying Parseval's theorem and then Taylor's theorem again, we arrive at the asymptotic
\begin{align*}
	W_{m, n}^+(x) &= \frac{T M}{\sqrt{\abs{x}}} e\Bigl(-\frac{x}{2 \pi} + \frac{T^2}{\pi x} \Bigr) \sum_{0 \leq l \leq L_1} \sum_{0 \leq l_1 \leq 2 l} \sum_{\frac{l_1}{4} \leq l_2 \leq L_2} c_{l,l_1,l_2} \frac{M^{2 l - l_1} T^{4 l_2 - l_1}}{x^{l + 3 l_2 - l_1}} \times {} \\
	&\qquad {} \times \widehat{k^*}^{(2 l - l_1)}\Bigl( - \frac{2 M T}{\pi x}\Bigr)  + O \biggl( \frac{T M}{\sqrt{\abs{x}}} \Bigl( \frac{T^4}{\abs{x}^3}\Bigr)^{L_2 + 1} + T \Bigl( \frac{M}{\sqrt{\abs{x}}}\Bigr)^{2 L_1 + 3} \biggr).
\end{align*}
Here $c_{l, l_1, l_2}$ are constants depending only on $l$, $l_1$, and $l_2$, where in particular from the above calculations $c_{0,0,0} = \frac{1 + i}{\sqrt{2}}$.

Now we are essentially done, because recalling how $\tilde H_{m, n}^{+, 1}(x)$ in \eqref{Hmndecomp} is composed of $W_{m, n}^+$, a similar term, and a remainder of $O(\frac{T \abs{x}}{M^8})$, we simply take $L_1$ and $L_2$ to be sufficiently large that the error $O(\frac{T \abs{x}}{M^8})$ dominates the error terms in $L_1$ and $L_2$.

It then suffices to study only the leading term $l = l_1 = l_2 = 0$, because the rest are of identical form and can be handled similarly.

From these asymptotics we infer that to study $\mathcal{R}^+_3$ it suffices to study
\begin{align*}
	 \tilde{\mathcal{R}}^+_3 &\coloneqq \frac{(1 + i) T M}{2 \pi p^{\frac{1}{4}}} \sum_{m \geq 1} \sum_{n \geq 1} \frac{A(n, m)}{m n^{\frac{3}{4}}} g\Bigl( \frac{m^2 n}{N} \Bigr) \times {} \\
	&\qquad {} \times  \sum_{c \leq \frac{C_2}{m}} \frac{S(n, p; c)}{c^{\frac{1}{2}}} e\Bigl( \frac{2 \sqrt{n p}}{c} - \frac{T^2 c}{4 \pi^2 \sqrt{n p}}\Bigr) \widehat{k^*} \Bigl( \frac{M T c}{2 \pi^2 \sqrt{n p}} \Bigr).\numberthis\label{R3tilde}
\end{align*}

The error in this asymptotic expansion is therefore of size
\begin{align*}
	 &\sum_{m \geq 1} \sum_{n \geq 1} \frac{A(n, m)}{m n^{\frac{1}{2}}} g \Bigl( \frac{m^2 n}{N} \Bigr) \sum_{c \leq \frac{C_2}{m}} c^{-1} S(n, p; c) \frac{T \abs{x}}{M^8} \\
	&\qquad\qquad\ll \sum_{m \geq 1} \sum_{n \geq 1} \frac{\abs{A(n, m)}}{m n^{\frac{1}{2}}} g \Bigl( \frac{m^2 n}{N} \Bigr) \sum_{c \leq \frac{C_2}{m}} c^{-\frac{3}{2} + \varepsilon} \frac{T n^{\frac{1}{2}} p^{\frac{1}{2}}}{M^8} \\
	&\qquad\qquad\ll T M^{-8}  p^{\frac{1}{2}} \sum_{m \geq 1} m^{-2} \sum_{n \geq 1} \abs{A(n, m)} g \Bigl( \frac{m^2 n}{N} \Bigr) C_2 \\
	&\qquad\qquad\ll T M^{-8} N p^{\frac{1}{2}} \frac{N^{\frac{1}{2}} p^{\frac{1}{2}}}{T^{1 - \varepsilon} M} = T^{\varepsilon} M^{-9} N^{\frac{3}{2}} p.\numberthis\label{asymperror}
\end{align*}
where we have used the Weil bound on the Kloosterman sum and recalling that $x = \frac{4 \pi \sqrt{n p}}{c}$ and $C_2 = \frac{N^{\frac{1}{2}} p^{\frac{1}{2}}}{T^{1 - \varepsilon} M}$.
Since the main term from $\mathcal{D}$ is of size $T M$, we need this to be at most $T^{1 - \varepsilon} M$.
This tells us how long of a sum $T^a < M \leq T^{1 - \varepsilon}$ we are allowed, namely, to use this asymptotic expansion, we need $a > \frac{7}{20} + \varepsilon$.
We will have occasion to tune this $a$ more finely later on, when we analyze $\tilde{\mathcal{R}}^+_3$ itself.

To see that we need to open the Kloosterman sum, notice how if we were to sum trivially over $n$, using Weil's bound of the Kloosterman sum, we get
\[
	\tilde{\mathcal{R}}^+_3 \ll T M C_2^{1 + \varepsilon} N^{\frac{1}{4}} \ll N^{\frac{3}{4}} \leq T^{\frac{9}{4} + \varepsilon},
\]
recalling how $N \leq T^{3 + \varepsilon}$.
Hence to save $T^{\frac{5}{4}} M^{-1}$ we open the Kloosterman sum and sum nontrivially over $n$ using the Voronoi formula for $\GL(3)$.

To this end we identify the terms depending on $n$ from the expression \eqref{R3tilde} for $\tilde{\mathcal{R}}_3^+$, except for the Fourier coefficients $A(n, m)$ and the Kloosterman sums $S(n, p; c)$ (since those will be handled by the Voronoi formula), as
\[
	\psi(y) = y^{-\frac{3}{4}} g\Bigl( \frac{m^2 y}{N} \Bigr) e \Bigl( \frac{2 \sqrt{y p}}{c} - \frac{T^2 c}{4 \pi^2 \sqrt{y p}} \Bigr) \widehat{k^*}\Bigl( \frac{M T c}{2 \pi^2 \sqrt{y p}} \Bigr).
\]

Opening the Kloosterman sum we get
\[
	\numberthis{\label{openkloosterman}}\sum_{n \geq 1} A(n, m) S(n, p; c) \psi(n) = \sum_{d \bar{d} \equiv 1 \pmod*{c}} e\Bigl( \frac{p \bar d}{c} \Bigr) \sum_{n \geq 1} A(n, m) e\Bigl( \frac{n d}{c} \Bigr) \psi(n).
\]
This is a productive thing to do because the Voronoi formula \eqref{voronoiformula} tells us precisely how to handle the inner sum, namely
\[
	\sum_{n \geq 1} A(n, m) e\Bigl( \frac{n d}{c} \Bigr) \psi(n) = c \sum_{\pm} \sum_{n_1 \mid c m} \sum_{n_2 \geq 1} \frac{A(n_1, n_2)}{n_1 n_2} S(m \bar d, \pm n_2; m c n_1^{-1}) \Psi^{\pm} \Bigl( \frac{n_2 n_1^2}{c^3 m} \Bigr).
\]

Recalling now $c \leq \frac{C_2}{m} = \frac{\sqrt{N p}}{T^{1 - \varepsilon} M m}$ and how $N \leq T^{3 + \varepsilon}$, the size of the argument after transforming times the length of the original $n$-sum is at least
\[
	\frac{n_2 n_1^2}{c^3 m} \frac{N}{m^2} \gg T^{\frac{21}{8} - \varepsilon} p^{-\frac{3}{2}} \gg 1
\]
since $p \ll T^{1 - \varepsilon}$. 
Consequently the version of \cite[Lemma~2.1]{li2011} in \eqref{voronoiformulaasymptotics} applies and the size of the integral transform $\Psi^\pm(x)$ reduces to studying
\begin{align}\label{lem2.1expansion}
	x^{\frac{2}{3}} c_1^{\pm} \int_0^\infty e(u_1(y)) a(y) \, d y + x^{\frac{2}{3}} d_1^{\pm} \int_0^\infty e(u_2(y)) a(y) \, d y,
\end{align}
where $c_1^\pm$ and $d_1^\pm$ are absolute constants depending on the Langlands parameters of $f$,
\[
	u_1(y) = \frac{2 \sqrt{y p}}{c} + 3 x^{\frac{1}{3}} y^{\frac{1}{3}}, \qquad u_2(y) = \frac{2 \sqrt{y p}}{c} - 3 x^{\frac{1}{3}} y^{\frac{1}{3}},
\]
and
\[
	a(y) = \psi(y) y^{-\frac{1}{3}} = g \Bigl( \frac{m^2 y}{N} \Bigr) \widehat{k^*}\Bigl( \frac{M T c}{2 \pi^2 \sqrt{y p}} \Bigr) e \Bigl( - \frac{T^2 c}{4 \pi^2 \sqrt{y p}} \Bigr) y^{-\frac{13}{12}}.
\]
As mentioned in the discussion of \eqref{voronoiformulaasymptotics}, there are more terms in $\Psi^\pm(x)$, but like in all our previous analysis they are of lower order and do not contribute anything of interest.

These are oscillatory integrals with a weight, so they lend themselves to stationary phase analysis.
In particular, the first integral with phase $u_1(y)$ has negligible contribution, since $u_1'(y) \gg c^{-1} y^{-\frac{1}{2}}$ and $a'(y) \ll T^2 c y^{-\frac{31}{12}}$ meaning that
\[
	u_1'(y) a'(y)^{-1} \gg c^{-2} T^{-2} y^{\frac{25}{12}} \gg T^{-\varepsilon} M^2 m^2 N^{-1} y^{\frac{25}{12}} \gg T^{2 a - \varepsilon}.
\]
Whatever we tune $a$ to be, this power of $T$ is no doubt positive, and so using partial integration as many times as we like we can make the contribution of this first integral bounded by as large a negative power of $T$ as we please.

The second integral is the one that contributes a stationary phase, and consequently the one we have to exercise more care with.
This time, because of the minus sign in the phase, we have
\[
	u_2'(y) = p^{\frac{1}{2}} c^{-1} y^{-\frac{1}{2}} - x^{\frac{1}{3}} y^{-\frac{2}{3}},
\]
meaning that if $x$ is bounded away from $p^{\frac{3}{2}} c^{-3} y^{\frac{1}{2}}$, in particular $x$ away from $p^{\frac{3}{2}} c^{-3} m^{-1} N^{\frac{1}{2}}$ since $y = n \sim N m^{-2}$, then we have no stationary phase, and the second integral has negligible contribution for the same reason as the first one.

Consequently to have a stationary phase we are looking at $x$ around the size $p^{\frac{3}{2}} c^{-3} m^{-1} N^{\frac{1}{2}}$, or in other words $n_2 \sim N^{\frac{1}{2}} p^{\frac{3}{2}} n_1^{-2}$.
In this range the relevant integral
\[
	\int_0^\infty e(u_2(y)) a(y) \, d y
\]
therefore has a stationary phase point at $y_0 = p^{-3} c^6 x^2$, which using the stationary phase method (\cite[Lemma~5.5.6]{huxley1996}) gives us
\begin{equation}\label{eq:stationaryphaseR3}
	\int_0^\infty e(u_2(y)) a(y) \, d y = \frac{a(y_0) e\Bigl(-x c^2 p^{-1} + \dfrac{1}{8}\Bigr)}{\sqrt{u_2''(y_0)}} + O( c^{\frac{7}{2}} T^4 N^{-\frac{11}{6}} m^{\frac{11}{3}} p^{-\frac{7}{4}} ).
\end{equation}
Note for ease of reference later that
\[
	u_2''(y) = \frac{2}{3} x^{\frac{1}{3}} y^{-\frac{5}{3}} - \frac{1}{2} p^{\frac{1}{2}} c^{-1} y^{-\frac{3}{2}},
\]
so in particular at $y_0 = p^{-3} c^6 x^2$ we have $u_2''(y_0) = \frac{1}{6} x^{-3} c^{-10} p^5$, so that
\[
	\frac{1}{\sqrt{u_2''(y_0)}} = \sqrt{6} x^{\frac{3}{2}} c^5 p^{- \frac{5}{2}} = \sqrt{6} n_2^{\frac{3}{2}} n_1^{3} c^{\frac{1}{2}} m^{-\frac{3}{2}} p^{-\frac{5}{2}}
\]
for $x = \frac{n_2 n_1^2}{c^3 m}$.

In order to now use this to control the error we need to do something about the Kloosterman sums that remain after using the Voronoi formula.
Per \eqref{openkloosterman} and the Voronoi formula we are left with
\[
	\sum_{d \bar d \equiv 1 \pmod*{c}} e\Bigl( \frac{p \bar d}{c} \Bigr) S(m \bar d, \pm n_2; m c n_1^{-1}).
\]
To handle this, open the Kloosterman sum and switch the order of summation to get
\begin{align*}
	\sum_{u \bar u \equiv 1 \pmod*{m c n_1^{-1}}} &e\Bigl( \frac{\pm n_2 \bar u}{m c n_1^{-1}} \Bigr) \sum_{d \bar d \equiv 1 \pmod*{c}} e\Bigl ( \frac{p \bar d}{c} \Bigr) e \Bigl( \frac{m \bar d u}{m c n_1^{-1}} \Bigr) \\
	&= 	\sum_{u \bar u \equiv 1 \pmod*{m c n_1^{-1}}} e\Bigl( \frac{\pm n_2 \bar u}{m c n_1^{-1}} \Bigr) \sum_{d \bar d \equiv 1 \pmod*{c}} e\Bigl( \frac{\bar d (p + u n_1)}{c} \Bigr) \\
	&= \sum_{u \bar u \equiv 1 \pmod*{m c n_1^{-1}}} e\Bigl( \frac{\pm n_2 \bar u}{m c n_1^{-1}} \Bigr) S(0, p + u n_1; c).
\end{align*}
The Kloosterman sum in the right-hand side,
\[
	S(0, a; c) = \sum_{v \bar v \equiv 1 \pmod*{c}} e\Bigl( \frac{a v}{c} \Bigr),
\]
is really just the Ramanujan sum, which is bounded by $(a, c)$, so the twisted sum of Kloosterman sums above is bounded by $m c^{1 + \varepsilon}$ since $n_1 \mid m c$.

This lets us bound the contribution to \eqref{R3tilde} from the error term in the stationary phase analysis by
\begin{align*}
	&M T p^{-\frac{1}{4}} \sum_{m \geq 1} m^{-1} \sum_{c \leq \frac{C_2}{m}} c^{\frac{1}{2}} \sum_{n_1 \mid c m} \sum_{n_2 \sim N^{\frac{1}{2}} p^{\frac{3}{2}} n_1^{-2}} \frac{\abs{A(n_1, n_2)}}{n_1 n_2}  \Bigl( \frac{n_2 n_1^2}{c^3 m} \Bigr)^{\frac{2}{3}} m c^{1 + \varepsilon} c^{\frac{7}{2}} T^4 N^{-\frac{11}{6}} m^{\frac{11}{3}} p^{-\frac{7}{4}} \\
	&\qquad\qquad\qquad \ll M^{-3} T^{1 + \varepsilon} N^{\frac{1}{2}} p^{1 + \varepsilon}
\end{align*}
using calculations similar to those of the error from the asymptotic expansion in \eqref{asymperror}.
Remember for this calculation how $x = \frac{n_2 n_1^2}{c^3 m}$ and how by \eqref{lem2.1expansion} there is an $x^{\frac{2}{3}}$ in front of the integral we are estimating with stationary phase analysis.

Since we want this to be dominated by $T M$, we again need to tune the length of $T^a < M \leq T^{1 - \varepsilon}$.
In particular, for $M^{-3} T^{1 + \varepsilon} N^{\frac{1}{2}} \leq M T^{-4a + 1 + \frac{3}{2} + \varepsilon}$ to be dominated by $T M$ we need $a > \frac{3}{8} + \varepsilon$.
Comparing this to $a > \frac{7}{20} + \varepsilon$ in order to get an admissible error from the asymptotic expansion before, we see that we have narrowed the range of $M$ slightly.

This means that, at this point, after using the Voronoi formula once and applying stationary phase analysis, we have
\begin{align*}
	 \tilde{\mathcal{R}}^+_3 &= \frac{T M p^{\frac{3}{4}}}{\pi} \sum_{m \geq 1} m^{-1} \sum_{c \leq \frac{C_2}{m}} c^{-1} \sum_{\pm} \sum_{n_1 \mid cm} n_1^{-1} \sum_{n_2 \geq 1} A(n_1, n_2) \times {} \\
	&\qquad {}\times \sum_{u \bar u \equiv 1 \pmod*{m c n_1^{-1}}} S(0, p + u n_1; c) e \Bigl( \frac{\pm n_2 \bar{u}}{m c n_1^{-1}} \Bigr) e \Bigl( \frac{-n_2 n_1^2}{c m p} \Bigr) b(n_2) \\
	&\qquad {} + O(M T^{\frac{5}{2} - 4 a + \varepsilon} p^{1 + \varepsilon})\numberthis{\label{afterVoronoi1}}
\end{align*}
where now
\[
	b(y) = y^{-1} g \Bigl( \frac{y^2 n_1^4}{N p} \Bigr) \widehat{k^*} \Bigl( \frac{ M T c m p }{2 \pi^2 y n_1^2} \Bigr) e \Bigl( \frac{-T^2 c m p}{4 \pi^2 y n_1^2} \Bigr) .
\]

If at this point we sum trivially over $n_2$, keeping in mind that $1 \leq \frac{n_2^2 n_1^4}{N p} \leq 2$, we find that the main term contributes about $T M C_2 \ll T M T^{\frac{1}{2}} M^{-1}$, so we need to save an additional $T^{\frac{1}{2} + \varepsilon} M^{-1}$ in order for this to be admissible.
To do this we have to once again switch the order of summation and use the Voronoi formula for a second time.

The relevant $n_2$-sum is
\[
	\sum_{n_2 \geq 1} A(n_1, n_2) e\Bigl( \frac{\pm n_2 \bar u}{m c n_1^{-1}} \Bigr) e\Bigl( \frac{-n_2 n_1^2}{c m p} \Bigr) b(n_2).
\]
We rewrite the character as
\begin{align*}
	e\Bigl( \frac{\pm n_2 \bar u}{m c n_1^{-1}} \Bigr) e\Bigl( \frac{-n_2 n_1^2}{c m p} \Bigr) &= e \Bigl( \frac{\pm n_2 \bar u}{m c n_1^{-1}} - \frac{n_2 n_1^2}{c m p} \Bigr) = e \Bigl( \frac{\pm n_2 \bar u p - n_2 n_1}{m c n_1^{-1} p} \Bigr) \\
	&= e\Bigl( \frac{n_2 (\pm \bar u p - n_1)}{m c n_1^{-1} p} \Bigr).
\end{align*}
Relabeling
\[
	\frac{\pm \bar u p - n_1}{m c n_1^{-1} p} = \frac{\bar u'}{c'}
\] 
with $(\bar u', c') = 1$ and $c' \mid m c n_1^{-1} p$, this gives us the following $n_2$-sum, expanded using the Voronoi formula:
\begin{align*}
	&\sum_{n_2 \geq 1} A(n_1, n_2) e\Bigl( \frac{n_2 \bar u'}{c'} \Bigr) b(n_2) \\
	&\qquad = c' \sum_{\pm} \sum_{l_1 \mid c' n_1} \sum_{l_2 \geq 1} \frac{A(l_1, l_2)}{l_1 l_2} S(n_1 u', \pm l_2; n_1 c' l_1^{-1}) B^{\pm}\Bigl( \frac{l_2 l_1^2}{c'^3 n_1} \Bigr). \numberthis{\label{Voronoi2}}
\end{align*}
Here $B^\pm(x)$ are the same integral transforms of $b(y)$ that $\Psi^\pm(x)$ are of $\psi(y)$ in our first use of the Voronoi formula, and so for the same reason we look at the size of the argument $x$ compared to the length of the sum prior to using the Voronoi formula in order to use \eqref{voronoiformulaasymptotics}:
\[
	\frac{l_2 l_1^2}{c'^3 n_1} \frac{N^{\frac{1}{2}} p^{\frac{1}{2}}}{n_1^2} \gg M^3 \gg 1
\]
and so again the asymptotics apply, and we can write the integral transform $B^\pm(x)$ in terms of
\begin{align*}
	B^\pm(x) = x^{\frac{2}{3}} c_1^\pm \int_0^\infty e(v_1(y)) q(y) \, d y - x^{\frac{2}{3}} d_1^\pm \int_0^\infty e(v_2(y)) q(y) \, d y \\ {} + \text{lower order terms}. \numberthis{\label{Bdecompose}}
\end{align*}
As before, $c_1^\pm$ and $d_1^\pm$ are absolute constants depending on $f$.
The phases are
\[
	v_1(y) = 3 x^{\frac{1}{3}} y^{\frac{1}{3}} - \frac{T^2 c m p}{4 \pi^2 y n_1^2} \quad \text{and} \quad v_2(y) = -3 x^{\frac{1}{3}} y^{\frac{1}{3}} - \frac{T^2 c m p}{4 \pi^2 y n_1^2},
\]
and
\[
	q(y) = y^{-\frac{4}{3}} g\Bigl( \frac{y^2 n_1^4}{N p} \Bigr) \widehat{k^*}\Bigl( \frac{M T c m p}{2 \pi^2 y n_1^2} \Bigr).
\]
We perform the same kind of analysis as we did after the first application of the Voronoi formula.
First,
\[
	v_1'(y) = x^{\frac{1}{3}} y^{-\frac{2}{3}} + \frac{T^2 c m p}{4 \pi^2 y^2 n_1^2} \gg \frac{T^2 c m p}{y^2 n_1^2}
\]
and $q'(y) \ll y^{-\frac{1}{3}} T^\varepsilon N^{-1} n_1^4 p^{-1}$, and hence
\[
	v_1'(y) q'(y)^{-1} \gg \frac{T^2 c m p}{y^2 n_1^2} y^{-\frac{1}{3}} T^\varepsilon N^{-1} n_1^4 p^{-1} \gg T^{2 - \varepsilon} c m p^{\frac{1}{2}} N^{-\frac{1}{2}} y^{\frac{4}{3}} \gg T^{\frac{1}{2} - \varepsilon}
\]
since $1 \leq \frac{y^2 n_1^4}{N p} \leq 2$.
Consequently, since that is a positive power of $T$, we can integrate by parts as many times as we like to get a sufficiently large negative power of $T$, so the integral with $v_1(y)$ as its phase is negligible.

Concerning the integral with phase $v_2(y)$, we have
\[
	v_2'(y) = -x^{\frac{1}{3}} y^{-\frac{2}{3}} + \frac{T^2 c m p}{4 \pi^2 y^2 n_1^2},
\]
so in order to have the same negligible asymptotics as in the previous case we need $\abs{v_2'(y)} \gg \frac{T^2 c m p}{y^2 n_1^2}$, or in other words we need $x$ away from
\[
	\frac{T^6 c^3 m^3 p n_1^2}{64 \pi^6 N^2}
\]
since $y$ is about $\frac{N^{\frac{1}{2}} p^{\frac{1}{2}}}{n_1^2}$.
Therefore if
\[
	x \geq \frac{T^6 c^3 m^3 p n_1^2}{10 \pi^6 N^2} \quad \text{or} \quad x \leq \frac{T^6 c^3 m^3 p n_1^2}{100 \pi^6 N^2}
\]
we have, by exactly the same argument as for $v_1(y)$, no stationary or small phase and the integral with phase $v_2(y)$ is negligible in this range.

It remains to study the integral with phase $v_2(y)$ in the range
\[
	\frac{T^6 c^3 m^3 p n_1^2}{100 \pi^6 N^2} \leq x \leq \frac{T^6 c^3 m^3 p n_1^2}{10 \pi^6 N^2}
\]
or equivalently the $l_2$-sum in the range
\[
	\frac{L_2}{100} \leq l_2 \leq \frac{L_2}{10}
\]
with
\[
	L_2 = \frac{T^6 c^3 m^3 p n_1^3 c'^3}{\pi^6 N^2 l_1^2}.
\]
Here we have
\[
	\abs{v_2''(y)} \gg \frac{T^2 c m p}{y^3 n_1^2} \gg T^2 c m p^{-\frac{1}{2}} N^{-\frac{3}{2}} n_1^4.
\]
By the Second derivative test (\cite[Lemma~5.1.3]{huxley1996}), we deduce from this and \eqref{Bdecompose} that
\begin{align*}
	B^\pm(x) &\ll x^{\frac{2}{3}} (T^2 c m p^{-\frac{1}{2}} N^{-\frac{3}{2}} n_1^4)^{-\frac{1}{2}} \Bigl( \frac{N^{\frac{1}{2}} p^{\frac{1}{2}}}{n_1^2} \Bigr)^{-\frac{4}{3}} T^{\varepsilon} \\
	&\ll T^{3+\varepsilon} c^{\frac{3}{2}} m^{\frac{3}{2}} p^{\frac{1}{4}} n_1^2 N^{-\frac{5}{4}}\numberthis{\label{B0asymptotics}}
\end{align*}
since at this point $x \ll T^6 c^3 m^3 p n_1^2 N^{-2}$.

Combining \eqref{afterVoronoi1}, \eqref{Voronoi2}, and \eqref{B0asymptotics}, along with the Weil bound $(n_2 c' l_1^{-1})^{\frac{1}{2} + \varepsilon}$ for the innermost Kloosterman sum and the bound $(p + u n_1, c)$ for the Ramanujan sum, we finally have
\begin{align*}
	\tilde{\mathcal{R}}^+_3 &\ll M T p^{\frac{3}{4}} \sum_{m \geq 1} m^{-1} \sum_{c \leq \frac{C_2}{m}} c^{-1} \sum_{n_1 \mid c m} n_1^{-1} \sum_{u \bar u \equiv 1 \pmod*{m c n_1^{-1}}} (p + u n_1, c) \\
	&\qquad {} \times c' \sum_{l_1 \mid c' n_1} \sum_{\frac{L_2}{100} \leq l_2 \leq \frac{L_2}{10}} \frac{\abs{A(l_1, l_2)}}{l_1 l_2} (n_1 c' l_1^{-1})^{\frac{1}{2} + \varepsilon} \\
	&\qquad {} \times T^{3+\varepsilon} c^{\frac{3}{2}} m^{\frac{3}{2}} p^{\frac{1}{4}} n_1^2 N^{-\frac{5}{4}} + O(M T^{\frac{5}{2} - 4 a + \varepsilon} p^{1 + \varepsilon}) \\
	&\ll M^{-3} T^{\varepsilon} N^{\frac{3}{4}} p + O(M T^{\frac{5}{2} - 4a + \varepsilon} p^{1 + \varepsilon}) \\
	&\ll M T^{\frac{9}{4} - 4 a + \varepsilon} p + O(M T^{\frac{5}{2} - 4a + \varepsilon} p^{1 + \varepsilon}).
\end{align*}
For this to be admissible, namely dominated by $T M$, we need $a > \frac{5}{16} + \varepsilon$ in  $T^a \leq M \leq T^{1 - \varepsilon}$, which of course is covered by $a > \frac{3}{8} + \varepsilon$ from earlier.

\section{Proof of \autoref{thm:maintheorem}: The $K$-Bessel function terms}\label{sec:Kbessel}

Again using partition of unity it suffices for $\mathcal{O}^-$ to consider
\[
	\mathcal{R}^- = \frac{1}{2} \sum_{m \geq 1} \sum_{n \geq 1} \frac{A(n, m)}{(m^2 n)^{\frac{1}{2}}} g\Bigl( \frac{m^2 n}{N} \Bigr) \sum_{c > 0} \frac{S(-n, p; c)}{c} H_{m, n}^- \Bigl ( \frac{4 \pi \sqrt{n p}}{c} \Bigr )
\]
where $g$ is a smooth function of compact support on $\interval{1}{2}$ and $N$ is at most $T^{3 + \varepsilon}$.

We prove the following:
\begin{lemma}
	\[
		\mathcal{R}^- = O(M^{-1} T^{\frac{3}{2} + \varepsilon} p^{\varepsilon}).
	\]
\end{lemma}

Recall how
\[
	H_{m, n}^-(x) = \frac{4}{\pi} \int_{-\infty}^\infty K_{2 i t}(x) \sinh(\pi t) k(t) V(m^2 n, t) t \, d t.
\]
As with $\mathcal{R}^+$ we will split the $c$-sum, but this time it suffices to split it into only two parts,
\[
	\mathcal{R}_1^- = \frac{1}{2} \sum_{m \geq 1} \sum_{n \geq 1} \frac{A(n, m)}{(m^2 n)^{\frac{1}{2}}} g\Bigl( \frac{m^2 n}{N} \Bigr) \sum_{c \geq \frac{C}{m}} \frac{S(-n, p; c)}{c} H_{m, n}^- \Bigl ( \frac{4 \pi \sqrt{n p}}{c} \Bigr )
\]
and
\[
	\mathcal{R}_2^- = \frac{1}{2} \sum_{m \geq 1} \sum_{n \geq 1} \frac{A(n, m)}{(m^2 n)^{\frac{1}{2}}} g\Bigl( \frac{m^2 n}{N} \Bigr) \sum_{c \leq \frac{C}{m}} \frac{S(-n, p; c)}{c} H_{m, n}^- \Bigl ( \frac{4 \pi \sqrt{n p}}{c} \Bigr ).
\]
As before we will tune the cut-off $C$ carefully in a moment.

The first part we deal with in essentially the same way we dealt with $\mathcal{R}_1^+$ involving the $J$-Bessel function, after first using an identity to switch from the $K$-Bessel function to the $I$-Bessel function (because we have similar integral representations for the $I$-Bessel function and the $J$-Bessel function).
That is to say, we move the line of integration in $H_{m, n}^-(x)$ carefully and use the integral representation for the Bessel function to extract bounds.

In particular, we have
\[
	K_\nu(z) = \frac{\pi}{2} \frac{I_{-\nu}(z) - I_\nu(z)}{\sin(\pi \nu)},
\]
so that
\begin{align*}
	H_{m, n}^-(x) &= 2 \int_{-\infty}^\infty \frac{I_{-2 i t}(x) - I_{2 i t}(x)}{\sin(2 i t \pi)} \sinh(\pi t) k(t) V(m^2 n, t) t \, d t \\
	&= -4 \int_{-\infty}^\infty \frac{I_{2 i t}(x)}{\sin(2 i t \pi)} \sinh(\pi t) k(t) V(m^2 n, t) t \, d t
\end{align*}
since $k(t)$ and $V(m^2 n, t)$ are even in $t$.

Moving the line of integration to $\Im t = -d$, which we will tune momentarily, this becomes
\begin{align*}
	-4 \int_{-\infty}^\infty \frac{I_{2d + 2 i y}(x)}{\sin((2 d + 2 i y) \pi)} \sinh(\pi(-d i + y)) k(-d i + y) &V(m^2 n, -d i + y) (-d i + y) \, d y +{} \\
	&{}+ \text{residue terms}.
\end{align*}
The residue terms come from denominator $\sin((2 u + 2 i y) \pi)$ with simple zeros at $2 u + 2 i y = l$, $l \in \Z$ and $l \neq 0$, and are of size
\[
	\ll I_l(x) V\Bigl(m^2 n, -i \frac{l}{2} \Bigr) \frac{l}{2}.
\]
For fixed $d$ there are only finitely many such residues for $l = 1, 2, \ldots < d$.

Per \cite[8.431~3]{gradshtein2015} we have the integral representation
\[
	I_\nu(z) = \frac{\Bigl(\dfrac{z}{2}\Bigr)^\nu}{\Gamma\Bigl(\nu + \dfrac{1}{2}\Bigr) \Gamma\Bigl(\dfrac{1}{2}\Bigr)} \int_0^\pi e^{\pm z \cos(\theta)} \sin(\theta)^{2 \nu} \, d \theta
\]
so that for $\Re(\nu) > -\frac{1}{2}$ we have the bound
\[
	I_\nu(z) \ll \Bigl ( \frac{z}{\abs{\Im \nu} + 1} \Bigr)^{\Re \nu} e^{\pi \Im \nu} e^z,
\]
and from \eqref{eq:StirlingV} we have $V(m^2 n, -d i + y) \ll \abs{y}^{3 d} (m^2 n)^{-d}$.
This makes the residue terms bounded by
\[
	\ll x^d (m^2 n)^{-\frac{1}{2}},
\]
and hence we have
\[
	H_{m, n}^-(x) \ll x^{2 d} e^x (m^2 n)^{-d} T^{d + 1 + \varepsilon} M + x^d (m^2 n)^{-\frac{1}{2}}
\]
since $y$ is about $T$ and the length of the integral is about $M$, both because of the exponential decay of $k(-d i + y)$.
This is the first step in tuning $C$, the lower bound on $c \geq \frac{C}{m}$.
In order for the $e^x$ factor to be insignificant we need $x = \frac{4 \pi \sqrt{n p}}{c}$ bounded above, say by $c \geq \sqrt{N}$.

With this we have, using the Weil bound on the Kloosterman sum,
\begin{align*}
	\mathcal{R}_1^- \ll \sum_{m \geq 1} \sum_{n \geq 1} &\frac{\abs{A(n, m)}}{(m^2 n)^{\frac{1}{2}}} g\Bigl( \frac{m^2 n}{N} \Bigr) \sum_{c \geq \frac{C}{m}} c^{-1} c^{\frac{1}{2} + \varepsilon} \times {} \\
	&{} \times \biggl(\Bigl( \frac{4 \pi \sqrt{n p}}{c} \Bigr)^{2 d} (m^2 n)^{-d} T^{d + 1 + \varepsilon} M + \Bigl( \frac{4 \pi \sqrt{n p}}{c} \Bigr)^{d} (m^2 n)^{-\frac{1}{2}} \biggr).
\end{align*}
The power of $c$ from the first term in the parentheses is $\frac{1}{2} - 1 - 2 d + \varepsilon$ and the power of $c$ from the second term in the parentheses is $\frac{1}{2} - 1 - d + \varepsilon$, so for the $c$-sum to converge we need $d > \frac{1}{2}$ when moving the line of integration.

This then makes the $c$-sum about $(\frac{C}{m})^{\frac{1}{2} - 2 d + \varepsilon}$, so by taking
\[
	C = \sqrt{N} + T,
\]
the $\sqrt{N}$ in case $m^2 \ll N$ is large and $T$ in case $m$ is small, and picking $d$ sufficiently large we make the entirety of $\mathcal{R}_1^-$ bounded by a large negative power of $T$.


Finally we handle
\[
	\mathcal{R}_2^- = \frac{1}{2} \sum_{m \geq 1} \sum_{n \geq 1} \frac{A(n, m)}{(m^2 n)^{\frac{1}{2}}} g\Bigl( \frac{m^2 n}{N} \Bigr) \sum_{c \leq \frac{C}{m}} \frac{S(-n, p; c)}{c} H_{m, n}^- \Bigl ( \frac{4 \pi \sqrt{n p}}{c} \Bigr )
\]
in a way exactly analogous to how $\mathcal{R}_2^+$ and $\mathcal{R}_3^+$ were handled.
We substitute the integral representation (from \cite[8.432~4]{gradshtein2015})
\[
	K_{2 i t}(x) = \frac{1}{2} \cosh(t \pi)^{-1} \int_{-\infty}^\infty \cos(x \sinh \zeta) e\Bigl( -\frac{t \zeta}{\pi} \Bigr) \, d\zeta
\]
into $H_{m, n}^-(x)$, make the change of variables $u = \frac{t - T}{M}$ we made with $\mathcal{R}_2^+$, extract the $e^{-\frac{(t + T)^2}{M^2}}$ part of $k(t)$ with negligible error because of exponential decay, identify the resulting integral in $u$ as the Fourier transformation of $k^*(u) = e^{-u^2} V(m^2 n, M u + T)$, make a second change of variables $\xi = -\frac{M \zeta}{\pi}$, and finally reduce all of it to studying the integral
\[
	H_{m, n}^{-,1}(x) = 4 T \int_{-\infty}^{\infty} \widehat{k^*}(\xi) \cos \Bigl( x \sinh \frac{\xi \pi}{M} \Bigr) e\Bigl( -\frac{T \xi}{M} \Bigr) \, d \xi.
\]
As with $\mathcal{R}_2^+$ there is a second term, the one arising from the $M u$ part of $t = M u + T$ instead of the $T$ part, but again like $\mathcal{R}_2^+$ that second term can be handled similarly.

The point of this is that we now have a oscillatory integral, so by looking at the phase
\[
	\phi(\xi) = -\frac{T \xi}{M} \pm \frac{x}{2 \pi} \sinh \frac{\xi \pi}{M}
\]
with derivative
\[
	\phi'(\xi) = -\frac{T}{M} \pm \frac{x}{2 M} \cosh \frac{\xi \pi}{M}
\]
we see that the integral can have a small or stationary phase if $\abs{x}$ is close to $T$, say $\frac{T}{100} \leq \abs{x} \leq 100 T$.
Hence for $x$ outside this range the integral is negligible, and we need only concern ourselves with $x$ inside this range.

The approach to get asymptotics for $H_{m, n}^{-, 1}(x)$ (again following \cite[Proposition~5.1]{li2011}) is the same as that for $H_{m, n}^{+, 1}(x)$, except slightly easier: we expand the $\sinh$ in the phase as a Taylor series, giving us
\begin{align*}
	\tilde H_{m, n}^{-,1}(x) &= 4 T \int_{-\infty}^\infty \widehat{k^*}(\xi) e\Bigl( -\frac{T \xi}{M} - \frac{x}{2 \pi} \sinh \frac{\xi \pi}{M} \Bigr) \, d \xi \\
	&= 4 T \int_{-\infty}^\infty \widehat{k^*}(\xi) e\Bigl( -\frac{T \xi}{M} + \frac{x \xi}{2 M} + \frac{\pi^2 x \xi^3}{12 M^3} + \frac{\pi^4 x \xi^5}{240 M^5} \Bigr) \, d \xi + O\Bigl( \frac{T \abs{x}}{M^7} \Bigr).
\end{align*}
for the plus part of the $\pm$ in the phase; the minus part is similar.

Expanding $e(\frac{\pi^2 x \xi^3}{12 M^3} + \frac{\pi^4 x \xi^5}{240 M^5} )$ as a Taylor series of order $L$ this gives us
\begin{align*}
	\tilde H_{m, n}^{-,1}(x) &= 4 T \int_{-\infty}^\infty \widehat{k^*}(\xi) e \Bigl( - \frac{(2 T - x) \xi}{2 M} \Bigr) \sum_{0 \leq l \leq L} \sum_{0 \leq l_1 \leq l} d_{l, l_1} \Bigl( \frac{x \xi^3}{M^3} \Bigr)^{l_1} \Bigl( \frac{x \xi^5}{M^5} \Bigr)^{l - l_1} \, d \xi + {} \\
	& \qquad\qquad\qquad\qquad{} + O\Bigl( \frac{T \abs{x}^{L + 1}}{M^{3 L + 3}} + \frac{T \abs{x}}{M^7} \Bigr)
\end{align*}
where $d_{l, l_1}$ are constants from the Taylor expansion.
In particular, $d_{0, 0} = 1$.
Identifying the integral in $\xi$ as a Fourier transform we get
\begin{align*}
	\tilde H_{m, n}^{-,1}(x) &= 4 T \sum_{0 \leq l \leq L} \sum_{0 \leq l_1 \leq l} d'_{l, l_1} \frac{x^l}{M^{5 l - 2 l_1}} {k^*}^{(5 l - 2 l_1)}\Bigl( \frac{x - 2 T}{2 M} \Bigr) + O\Bigl( \frac{T \abs{x}^{L + 1}}{M^{3 L + 3}} + \frac{T \abs{x}}{M^7} \Bigr)
\end{align*}
where $d'_{l, l_1} = d_{l, l_1} (2 \pi i)^{-5 l + 2 l_1}$ are relabeled constants.
Again in particular $d'_{0, 0} = 1$.

Hence, like with $\mathcal{R}_3^+$, we take $L$ to be sufficiently large that the second error term dominates the error, and study the $l = l_1 = 0$ term since the other terms can be treated similarly.

This means that to study $\mathcal{R}_2^-$ it suffices to study
\[
	\tilde{\mathcal{R}}_2^- \coloneqq T \sum_{m \geq 1} \sum_{n \geq 1} \frac{A(n, m)}{(m^2 n)^{\frac{1}{2}}} g\Bigl( \frac{m^2 n}{N} \Bigr) \sum_{\frac{\sqrt{N p}}{100 T m} \leq c \leq \frac{100 \sqrt{N p}}{T m}} \frac{S(-n, p; c)}{c} k^* \Biggl( \frac{\dfrac{4 \pi \sqrt{n p}}{c} - 2 T}{2 M} \Biggr).
\]

The error term from these asymptotics on $H_{m, n}^-(x)$ is of size $O(\frac{T \abs{x}}{M^7})$, so to $\mathcal{R}_2^-$ as a whole this contributes
\begin{align*}
	&\sum_{m \geq 1} \sum_{n \geq 1} \frac{A(n, m)}{m n^{\frac{1}{2}}} g\Bigl( \frac{m^2 n}{N} \Bigr) \sum_{\frac{\sqrt{N p}}{100 T m} \leq c \leq \frac{100 \sqrt{N p}}{T m}} \frac{S(-n, p; c)}{c} \frac{T \abs{x}}{M^7} \\
	&\qquad\qquad\ll \sum_{m \geq 1} \sum_{n \geq 1} \frac{\abs{A(n, m)}}{m n^{\frac{1}{2}}} g\Bigl( \frac{m^2 n}{N} \Bigr) \sum_{\frac{\sqrt{N p}}{100 T m} \leq c \leq \frac{100 \sqrt{N p}}{T m}} c^{-\frac{3}{2} + \varepsilon} \frac{T n^{\frac{1}{2}} p^{\frac{1}{2}}}{M^7} \\
	&\qquad\qquad\ll T M^{-7} p^{\frac{1}{2}} \sum_{m \geq 1} m^{-1} \sum_{n \geq 1} \abs{A(n, m)} g\Bigl( \frac{m^2 n}{N} \Bigr) \Bigl( \frac{\sqrt{N p}}{T m} \Bigr)^{-\frac{1}{2} + \varepsilon} \\
	&\qquad\qquad\ll T^{\frac{3}{2} - \varepsilon} M^{-7} N^{-\frac{1}{4} + \varepsilon} p^{\frac{1}{4} + \varepsilon} \sum_{m \geq 1} m^{-\frac{1}{2} - \varepsilon} m \frac{N}{m^2} = T^{\frac{3}{2} - \varepsilon} M^{-7} N^{\frac{3}{4} + \varepsilon} p^{\frac{1}{4} + \varepsilon}.
\end{align*}
For this to be admissible we need $a > \frac{11}{32} + \varepsilon$, which is satisfied by the requirement $a > \frac{3}{8} + \varepsilon$ already established from $\tilde{\mathcal{R}}_3^+$.

If at this point we sum $\tilde{\mathcal{R}}_2^-$ trivially over $n$ using the Weil bound for the Kloosterman sum, we get
\[
	\tilde{\mathcal{R}}_2^- \ll T^{\frac{1}{2}} N^{\frac{3}{4} + \varepsilon} \ll T^{\frac{11}{4} + \varepsilon}.
\]
Hence we must save $T^{\frac{7}{4} + \varepsilon} M^{-1}$, and in order to do so we must repeat what we did previously: open the Kloosterman sum, apply the Voronoi formula for $\GL(3)$, and analyze the size of the resulting integral transform.

The relevant $n$-sum in $\tilde{\mathcal{R}}_2^-$ is
\[
	\sum_{n \geq 1} A(n, m) S(-n, p; c) \psi(n) = \sum_{d \bar d \equiv 1 \pmod*{c}} e\Bigl( \frac{-p \bar d}{c} \Bigr) \sum_{n \geq 1} A(n, m) e\Bigl( \frac{n d}{c} \Bigr) \psi(n)
\]
where
\[
	\psi(y) = g\Bigl( \frac{m^2 y}{N} \Bigr) k^* \Biggl( \frac{\dfrac{4 \pi \sqrt{y p}}{c} -2 T}{2 M} \Biggr) y^{-\frac{1}{2}},
\]
and the Voronoi formula says the inner sum is
\begin{align*}
	&\sum_{n \geq 1} A(n, m) e\Bigl( \frac{n d}{c} \Bigr) \psi(n) \\
	&\qquad = c \sum_{\pm} \sum_{n_1 \mid c m} \sum_{n_2 \geq 1} \frac{A(n_1, n_2)}{n_1 n_2} S(m \bar d, \pm n_2; m c n_1^{-1}) \Psi^{\pm}\Bigl( \frac{n_2 n_1^2}{c^3 m}\Bigr)
\end{align*}
As before, we look at the size of the argument in the integral transforms times the length of the original $n$-sum:
\[
	\frac{n_2 n_1^2}{c^3 m} \frac{N}{m^2} \gg \frac{n_2 n_1^2}{m^3} \frac{N T^3 m^3}{N^{\frac{3}{2}}} \gg \frac{T^3}{N^{\frac{1}{2}} p^{\frac{3}{2}}} \gg T^{\varepsilon} \gg 1
\]
since now $c \leq \frac{100 N^{\frac{1}{2}} p^{\frac{1}{2}}}{T m}$ and $p \ll T^{1 - \varepsilon}$.
Therefore by \eqref{voronoiformulaasymptotics}, to study $\Psi^{\pm}(x)$ for $x = \frac{n_2 n_1^2}{c^3 m}$, it suffices to consider
\begin{equation}\label{eq:psiasymp}
	x^{\frac{2}{3}} \int_0^\infty e(\phi(y)) a(y) \, d y
\end{equation}
with the phase
\[
	\phi(y) = 3 x^{\frac{1}{3}} y^{\frac{1}{3}},
\]
and amplitude
\[
	a(y) = \psi(y) y^{-\frac{1}{3}} = g\Bigl( \frac{m^2 y}{N} \Bigr) k^* \Biggl( \frac{\dfrac{4 \pi \sqrt{y p}}{c} -2 T}{2 M} \Biggr) y^{-\frac{5}{6}}.
\]
The analysis this time is slightly easier than in the case for the $J$-Bessel function terms because this time $\psi(y)$ does not involve an exponential causing oscillations.

Since $\phi'(y) = x^{\frac{1}{3}} y^{-\frac{2}{3}}$ and $a'(y) \ll y^{-\frac{11}{6}}$, we get
\[
	\phi'(y) (a'(y))^{-1} \gg x^{\frac{1}{3}} y^{\frac{7}{6}} \gg n_2^{\frac{1}{3}} n_1^{\frac{2}{3}} c^{-1} m^{-3} N^{\frac{7}{12}}
\]
since $x = \frac{n_2 n_1^2}{c^3 m}$ and $y = n \sim N m^{-2}$ because of the factor of $g(\frac{m^2 y}{N})$ in $\psi(y)$.
Now if
\[
	n_2 \gg \frac{N^{\frac{1}{2}} T^{\varepsilon}}{M^3 n_1^2},
\]
this is makes
\[
	\phi'(y) (a'(y))^{-1} \gg T^\varepsilon,
\]
meaning that by integrating by parts many times, the contribution to $\tilde{\mathcal{R}}_2^-$ from such $n_2$ is negligible.
We therefore turn to
\[
	n_2 \ll \frac{N^{\frac{1}{2}} T^{\varepsilon}}{M^3 n_1^2}.
\]
Recalling how $k^*(u) = e^{-u^2} V(m^2 n, M u + T) \ll (1 + \abs{u})^{-A}$ for any $A > 0$, we have that $a(y)$ is negligible unless the argument
\[
	u = \frac{\dfrac{4 \pi \sqrt{y p}}{c} -2 T}{2 M}
\]
of $k^*(u)$ is about zero, i.e.,
\[
	\abs[\Bigg]{\frac{\dfrac{4 \pi \sqrt{y p}}{c} -2 T}{2 M}} = \abs[\Bigg]{\frac{\dfrac{2 \pi \sqrt{y p}}{c} - T}{M}} \leq T^\varepsilon.
\]
From this we deduce
\[
	\frac{1}{4 \pi^2 p} (T c - T^\varepsilon M c)^2 \leq y \leq 	\frac{1}{4 \pi^2 p} (T c + T^\varepsilon M c)^2.
\]
Therefore the length of the integral \eqref{eq:psiasymp} is about $T^{1 + \varepsilon} M c^2 p^{-1}$, the integrand is about $y^{-\frac{5}{6}}$, and there is an $x^{\frac{2}{3}}$ in front of the integral, so
\[
	\Psi^\pm(x) \ll x^{\frac{2}{3}} y^{-\frac{5}{6}} T^{1 + \varepsilon} M c^2 p^{-1} \ll x^{\frac{2}{3}} (N m^{-2})^{-\frac{5}{6}} T^{1 + \varepsilon} M c^2 p^{-1}.
\]
Plugging this into $\tilde{\mathcal{R}}_2^-$, along with the estimates on the Kloosterman sum coming from the Voronoi formula, we then get
\begin{align*}
	\tilde{\mathcal{R}}_2^- &\ll T \sum_{m \leq \sqrt{N}} m^{-1} \sum_{\frac{\sqrt{N p}}{100 T m} \leq c \leq \frac{100 \sqrt{N p}}{T m}} \sum_{n_1 \mid c m} \sum_{n_2 \ll \frac{\sqrt{N} T^{\varepsilon}}{M^3 n_1^2}} \frac{\abs{A(n_1, n_2)}}{n_1 n_2} m c^{1 + \varepsilon} \times {} \\
	& \qquad\qquad\qquad\qquad {} \times \Bigl( \frac{n_2 n_1^2}{c^3 m} \Bigr)^{\frac{2}{3}} (N m^{-2})^{-\frac{5}{6}} T^{1 + \varepsilon} M c^2 p^{-1} \\
	&= T^{2 + \varepsilon} M N^{-\frac{5}{6}} p^{-1} \sum_{m \leq \sqrt{N}} m \sum_{\frac{\sqrt{N p}}{100 T m} \leq c \leq \frac{100 \sqrt{N p}}{T m}} c^{1 + \varepsilon} \sum_{n_1 \mid c m} n_1^{\frac{1}{3}} \sum_{n_2 \ll \frac{\sqrt{N} T^{\varepsilon}}{M^3 n_1^2}} n_2^{-\frac{1}{3}} \abs{A(n_1, n_2)} \\
	&\ll T^{2 + \varepsilon} M N^{-\frac{5}{6}} p^{-1} \sum_{m \leq \sqrt{N}} m \sum_{\frac{\sqrt{N p}}{100 T m} \leq c \leq \frac{100 \sqrt{N p}}{T m}} c^{1 + \varepsilon} \sum_{n_1 \mid c m} n_1^{\frac{1}{3}} \Bigl( \frac{\sqrt{N} T^{\varepsilon}}{M^3 n_1^2} \Bigr)^{\frac{2}{3}} n_1 \\
	&\ll T^{2 + \varepsilon} M^{-1} N^{-\frac{1}{2}} p^{-1} \sum_{m \leq \sqrt{N}} m^{1 + \varepsilon} \sum_{\frac{\sqrt{N p}}{100 T m} \leq c \leq \frac{100 \sqrt{N p}}{T m}} c^{1 + \varepsilon} \\
	&\ll T^{\varepsilon} N^{\frac{1}{2} + \varepsilon} M^{-1} p^{\varepsilon} \ll M T^{\frac{3}{2} - 2 a + \varepsilon} p^{\varepsilon},
\end{align*}
which is admissible if $a > \frac{1}{4}$, and so covered by $a > \frac{3}{8}$.

This finishes the calculation for $\tilde{\mathcal{R}}_2^-$, and hence for $\mathcal{R}^-$ as a whole, and therefore of the entire moment.
Hence this finishes the proof of \autoref{thm:maintheorem}.

\section{Proof of \autoref{thm:derivativemoment}}\label{sec:derivativemoment}

The calculations for \autoref{thm:derivativemoment}, the twisted first moment of the derivative $L'(\frac{1}{2}, f \times u_j)$ for $\Set{u_j}$ an orthonormal basis of odd Hecke--Maass forms for $\SL(2, \Z)$, are essentially the same as the ones above.
The main difference is that of the diagonal main term, because this time the integral transform in the approximate functional equation has an order 2 pole instead of a simple pole.

For the derivative $L'(\frac{1}{2}, f \times u_j)$ we have the approximate functional equation
\begin{align*}
	L'\Bigl( \frac{1}{2}, f \times u_j \Bigr) &= \sum_{m \geq 1} \sum_{n \geq 1} \frac{\lambda_j(n) A(n, m)}{(m^2 n)^{\frac{1}{2}}} U_-(m^2 n, t_j) + {} \\
	&\qquad\qquad{} + \sum_{m \geq 1} \sum_{n \geq 1} \frac{\lambda_j(n) A(m, n)}{(m^2 n)^{\frac{1}{2}}} U_+(m^2 n, t_j)\numberthis{\label{eq:Lprimeapproxfunc}}
\end{align*}
where
\[
	U_{\mp}(y, t) = \frac{1}{2 \pi i} \int_{(1000)} y^{-u} F(u) \frac{\gamma_{\mp}\Bigl(\dfrac{1}{2} + u, t\Bigr)}{\gamma_-\Bigl(\dfrac{1}{2}, t\Bigr)} \frac{d u}{u^2}.
\]
As in the the proof \autoref{thm:maintheorem} we will focus on the $A(n, m)$ terms coming from the first sum, the dual sum one being treated similarly.
We consequently suppress the subscripts $\mp$ from now on.

We have the Kuznetsov trace formula for odd Maass forms (see \cite[Section~3]{conreyiwaniec2000}):
\begin{lemma}
	Let $h(t)$ be even, holomorphic and bounded $h(t) \ll (\abs{t} + 1)^{-2 - \varepsilon}$ in the strip $\abs{\Im t} \leq \frac{1}{2} + \varepsilon$.
	Then
	\begin{align*}
	  &\sideset{}{^*} \sum_j h(t_j) \omega_j \lambda_j(m) \lambda_j(n) \\
	  &= \frac{1}{2} \delta(m, n) H + \sum_{c > 0} \frac{1}{2 c} \left ( S(m, n; c) H^+\Bigl ( \frac{4 \pi \sqrt{m n}}{c} \Bigr ) - S(-m, n; c) H^- \Bigl ( \frac{4 \pi \sqrt{m n}}{c} \Bigr ) \right ),
	\end{align*}
	where $^*$ restricts the sum to odd Maass forms,
	\[
		\omega_j = \frac{4 \pi \abs{\rho_j(1)}^2}{\cosh(\pi t_j)},
	\]
	\[
		\omega(t) = \frac{4 \pi \abs[\Big]{\phi\Bigl(1, \dfrac{1}{2} + i t\Bigr)}^2}{\cosh(\pi t)}.
	\]
	\[
		H = \frac{1}{\pi} \int_{-\infty}^\infty h(t) \tanh(\pi t) t \, d t,
	\]
	\[
		H^+(x) = 2 i \int_{-\infty}^\infty J_{2 i t}(x) \frac{h(t) t}{\cosh(\pi t)} \, d t,
	\]
	and
	\[
		H^-(x) = \frac{4}{\pi} \int_{-\infty}^\infty K_{2 i t}(x) \sinh(\pi t) h(t) t \, d t.
	\]
\end{lemma}

Applying this and \eqref{eq:Lprimeapproxfunc} to the twisted moment in \autoref{thm:derivativemoment} we get
\begin{gather*}
	\sideset{}{^*}\sum_j k(t_j) \omega_j \lambda_j(p) L'\Bigl( \frac{1}{2}, f \times u_j \Bigr) = \frac{1}{2} \sum_{m \geq 1} \sum_{n \geq 1} \frac{A(n, m)}{(m^2 n)^{\frac{1}{2}}} \Biggl ( \delta(n, p) H_{m, n} + {} \\
	{} + \sum_{c > 0} \frac{1}{c} \biggl (  S(n, p; c) H_{m, n}^+\Bigl ( \frac{4 \pi \sqrt{n p}}{c} \Bigr ) - S(-n, p; c) H_{m, n}^- \Bigl ( \frac{4 \pi \sqrt{n p}}{c} \Bigr ) \biggr ) \Biggr ) + {} \\
	{} + \text{dual sum}
\end{gather*}
where this time
\[
  H_{m, n} = \frac{1}{\pi} \int_{-\infty}^\infty k(t) U(m^2 n, t) \tanh(\pi t) t \, d t,
\]
\[
	H_{m, n}^+(x) = 2 i \int_{-\infty}^\infty J_{2 i t}(x) \frac{k(t) U(m^2 n, t) t}{\cosh(\pi t)} \, d t,
\]
and
\[
	H_{m, n}^-(x) = \frac{4}{\pi} \int_{-\infty}^\infty K_{2 i t}(x) \sinh(\pi t) k(t) U(m^2 n, t) t \, d t.
\]

As before we split the resulting sum into diagonal and off-diagonal terms,
\[
	\mathcal{D} + \mathcal{O}^+ + \mathcal{O}^-,
\]
where
\[
	\mathcal{D} = \frac{1}{2} \sum_{m \geq 1} \sum_{n \geq 1} \frac{A(n, m)}{(m^2 n)^{\frac{1}{2}}} \delta(n, p) H_{m, n},
\]
\[
	\mathcal{O^+} = \frac{1}{2} \sum_{m \geq 1} \sum_{n \geq 1} \frac{A(n, m)}{(m^2 n)^{\frac{1}{2}}} \sum_{c > 0} \frac{S(n, p; c)}{c} H_{m, n}^+ \Bigl (\frac{4 \pi \sqrt{n p}}{c} \Bigr ),
\]
and
\[
	\mathcal{O}^- = \frac{1}{2} \sum_{m \geq 1} \sum_{n \geq 1} \frac{A(n, m)}{(m^2 n)^{\frac{1}{2}}} \sum_{c > 0} \frac{S(-n, p; c)}{c} H_{m, n}^- \Bigl ( \frac{4 \pi \sqrt{n p}}{c} \Bigr ).
\]

Note that $U(y, t)$ is bounded by $\log(t^3 y^{-1})$.
Since this is a power savings compared to the bound we used for $V(m^2 n, -\frac{1}{2} + y)$ before, we see that exactly the same calculations for the off-diagonal terms will work out for the derivative, and so the off-diagonal terms here contribute precisely the same error terms as they did in \autoref{thm:maintheorem}.

We hence turn our focus toward the diagonal term.
By the same shuffling of Fourier coefficients using the Hecke relation as in \eqref{diagonalterms}, we first have
\begin{align*}
	\mathcal{D} &= \frac{A(p, 1)}{2 p^{\frac{1}{2}}} \sum_{m \geq 1} \frac{A(1, m)}{m} H_{m, p} - \frac{1}{2 p^{\frac{3}{2}}} \sum_{m \geq 1} \frac{A(1, m)}{m} H_{m p, p}.
\end{align*}
Following the same approach as in Section~\ref{sec:diagonal}, we bring the sum into $H_{m, p}$ and get (for the first sum, the second one is handled similarly),
\begin{align*}
	\sum_{m \geq 1} \frac{A(1, m)}{m} H_{m, p} &= \frac{1}{\pi} \int_{-\infty}^\infty k(t) \Bigl( \sum_{m \geq 1} \frac{A(1, m)}{m} U(m^2 p, t) \Bigr) \tanh(\pi t) t \, d t.
\end{align*}
Note how the remaining exponential in the integral makes the length of the integral about $M$, and $t$ is about $T$.

Next, to estimate the inner sum we bring the sum all the way into the integral defining $U(m^2 p, t)$, i.e.,
\[
	\sum_{m \geq 1} \frac{A(1, m)}{m} U(m^2 p, t) = \frac{1}{2 \pi i} \int_{(1000)} p^{-u} L(1 + 2 u, \tilde f) F(u) \frac{\gamma\Bigl(\dfrac{1}{2} + u, t\Bigr)}{\gamma\Bigl(\dfrac{1}{2}, t\Bigr)} \frac{d u}{u^2}.
\]
Move the line of integration to $\Re(u) = -\frac{1}{7} + \varepsilon$, past the order 2 pole at $u = 0$ and we pick up the residue
\[
	2 L'(1, \tilde f) + 3 L(1, \tilde f) \log\abs{t} - 3 L(1, \tilde f) \log(2 \pi) - L(1, \tilde f) \log p + O(\abs{t}^{-1} p^{\varepsilon}).
\]
By Stirling's formula and the convexity bound for $L(s, \tilde f)$ we therefore get
\begin{align*}
	\sum_{m \geq 1} &\frac{A(1, m)}{m} U(m^2 p, t) \\
	&= 3 L(1, \tilde f) \log\abs{t} + 2 L'(1, \tilde f) - 3 L(1, \tilde f) \log(2 \pi) - L(1, \tilde f) \log p + {} \\
	&\qquad{} + O(\abs{t}^{-1} p^{\varepsilon}) + O(p^{\frac{1}{7} - \varepsilon} \abs{t}^{-\frac{6}{7} + \varepsilon} )
\end{align*}
for $t$ about $T$.

Plugging this back into the integral for $H_{m, p}$, we therefore get
\begin{align*}
	\sum_{m \geq 1} \frac{A(1, m)}{m} H_{m, p} &= \frac{3 L(1, \tilde f)}{\pi} \int_{-\infty}^{\infty} k(t) \tanh(\pi t) t \log\abs{t} \, d t + {}\\
	&\qquad{} + K \int_{-\infty}^\infty k(t) \tanh(\pi t) t \, d t + O(T^{\frac{1}{7} + \varepsilon} M p^{\frac{1}{7} - \varepsilon})
\end{align*}
where $K = 2 L'(1, \tilde f) - 3 L(1, \tilde f) \log(2 \pi) - L(1, \tilde f) \log p$.

The sum with $H_{m p, p}$ in place of $H_{m, p}$ works out precisely the same.

Hence
\begin{align*}
  \mathcal{D} &= \frac{3 L(1, \tilde f) \bigl( A(p, 1) p - 1\bigr)}{2 p^{\frac{3}{2}} \pi} \int_{-\infty}^{\infty} k(t) \tanh(\pi t) t \log\abs{t} \, d t + {} \\
	&\qquad {} + \frac{K \bigl( A(p, 1) p - 1\bigr)}{2 p^{\frac{3}{2}} \pi} \int_{-\infty}^{\infty} k(t) \tanh(\pi t) t \, d t + O(T^{\frac{1}{7} + \varepsilon} M p^{\varepsilon})
\end{align*}
with $K = 2 L'(1, \tilde f) - 3 L(1, \tilde f) \log(2 \pi) - L(1, \tilde f) \log p$.

The dual sum consequently works out identically, only with switching $\tilde f$ for $f$ and switching the Fourier coefficients, so
\begin{align*}
	\tilde{\mathcal{D}} &= \frac{3 L(1, f) \bigl( A(1, p) p - 1\bigr)}{2 p^{\frac{3}{2}} \pi} \int_{-\infty}^{\infty} k(t) \tanh(\pi t) t \log\abs{t} \, d t + {} \\
	&\qquad {} + \frac{\tilde K \bigl( A(p, 1) p - 1\bigr)}{2 p^{\frac{3}{2}} \pi} \int_{-\infty}^{\infty} k(t) \tanh(\pi t) t \, d t + O(T^{\frac{1}{7} + \varepsilon} M p^{\varepsilon})
\end{align*}
where $\tilde K = 2 L'(1, f) - 3 L(1, f) \log(2 \pi) - L(1, f) \log p$.

Together with the previously computed bounds for the off-diagonal terms which as discussed still apply, this proves \autoref{thm:derivativemoment}.

%
%
%
%

\section*{Acknowledgements}

The author would like to thank Sheng-Chi Liu for his tireless support and helpful insight, comments and encouragement, and the anonymous referee for their detailed and valuable feedback.

\bibliography{references}
\bibliographystyle{abbrv}

\end{document}